\documentclass[11pt]{article}
\usepackage{}
\usepackage{amsmath,amsfonts}
\usepackage{txfonts,wasysym,lscape,amssymb,mathrsfs,extarrows}
\usepackage[all,pdf]{xy}
\usepackage{graphicx}
\usepackage{subfigure}
 \usepackage{multirow}
 \usepackage{rotating}
 \usepackage[table,xcdraw]{xcolor}
\RequirePackage{color,graphicx}

\newtheorem{theorem}{Theorem}[section]
\newtheorem{proposition}{Proposition}[section]
\newtheorem{lemma}{Lemma}[section]

\newtheorem{algorithm}{Algorithm}[section]

\newtheorem{example}{Example}[section]
\newtheorem{remark}{Remark}[section]

\newtheorem{assumption}{Assumption}[section]

\title{Majorized Semi-proximal Alternating Coordinate Method for  Nonsmooth Convex-Concave Minimax Optimization}

\author{Yu-Hong Dai\footnote{LSEC, ICMSEC, AMSS, Chinese Academy of Sciences, Beijing 100190, P. R. China. {\sl Email}: dyh@lsec.cc.ac.cn.
This author was supported by the Natural Science Foundation of China (Nos. 11991020, 11631013, 12021001, 11971372 and 11991021) and
Beijing Academy of Artificial Intelligence (BAAI).}
 \footnote{School of Mathematical Sciences, University of Chinese Academy of Sciences, Beijing 100049, P. R. China.}\quad Jiani Wang\footnote{School of Mathematical Sciences, Dalian University of Technology, Dalian 116024, P. R. China. {\sl Email}: 1269338213@qq.com.} \quad and \quad Liwei Zhang \footnote{ School of Mathematical Sciences, Dalian University of Technology, Dalian 116024, P. R. China. {\sl Email}: lwzhang@dlut.edu.cn. This author was supported by the Natural Science Foundation of China (Nos. 11971089 and 11731013).}}

\date{}
\begin{document}

\maketitle

\begin{abstract}
Minimax optimization problems are an important class of optimization problems arising from modern machine learning and traditional research areas. While there have been many numerical algorithms for solving smooth convex-concave minimax problems, numerical algorithms for nonsmooth convex-concave minimax problems are very rare. This paper aims to develop an efficient numerical algorithm for a structured nonsmooth convex-concave minimax problem. A majorized semi-proximal alternating coordinate method (mspACM) is proposed, in which a majorized quadratic convex-concave function is adopted for approximating the smooth part of the objective function and semi-proximal terms are added in each subproblem. This construction enables the subproblems at each iteration are solvable and even easily solved when the semiproximal terms are cleverly chosen. We prove the global convergence of the algorithm mspACM under mild assumptions, without requiring strong convexity-concavity condition. Under the locally metrical subregularity of the solution mapping, we prove that the algorithm mspACM has the linear rate of convergence.  Preliminary numerical results are reported to verify the efficiency of the algorithm mspACM.
\vskip 6 true pt \noindent \textbf{Key words}: minimax optimization, convexity-concavity, semiproximal term, global convergence, rate of convergence, locally metrical subregularity.
\vskip 12 true pt \noindent \textbf{AMS subject classification}: 90C30
\end{abstract}
\bigskip\noindent

\section{Problem setting}
 \setcounter{equation}{0}
\quad \, Let ${\cal X}$ and ${\cal Y}$ be two finite-dimensional real Hilbert spaces equipped
with a scalar product $\langle \cdot,\cdot \rangle$ and its induced norm $\|\cdot\|$. Let $K: {\cal X} \times {\cal Y} \rightarrow \Re$ be a continuously differentiable convex-concave function, $f: {\cal X} \rightarrow \overline{\Re}$ and $g: {\cal Y} \rightarrow \overline{\Re}$ be proper lower semicontinuous convex functions.
We consider the following nonsmooth  minimax optimization problem
\begin{equation}\label{cminimax}
\min_{x \in {\cal X}}\max_{y \in {\cal Y}}L(x,y):=f(x)+K(x,y)-g(y).
\end{equation}
The mathematical model (\ref{cminimax}) covers a lot of interesting convex-concave minimax problems appeared in the literature. We only list two examples here.
\begin{example}\label{ex-1}
Let $X \subset {\cal X}$ and $Y \subset {\cal Y}$ be two closed convex sets. The following constrained minimax optimization problem are frequently studied in the literature.
\begin{equation}\label{Ecminimax}
\min_{x \in X}\max_{y \in Y}K(x,y).
\end{equation}
Obviously
the constrained minimax optimization problem
(\ref{Ecminimax}) can be written as the form of  Problem (\ref{cminimax}) if we set $f(x)=\delta_X(x)$ and $g(y)=\delta_Y(y)$.
\end{example}
\begin{example}\label{ex-2}
In Example (\ref{ex-1}), let $X \subset {\cal X}$ and $Y \subset {\cal Y}$ be two closed convex sets specified as
$$
X=\{x\in {\cal X}: G(x) \in C\},\quad Y=\{y\in {\cal Y}: H(y) \in D\},
$$
where $C$ and $D$ are closed convex sets in some finite dimensional spaces and the set-valued mappings
$$
x:\rightarrow G(x)-C \quad
\mbox{   and   } \quad
y:\rightarrow H(y)-D
$$
are graph-convex (under this assumption $X$ and $Y$ are convex sets). Then
the constrained minimax optimization problem
(\ref{Ecminimax}) can be written as the form of  Problem (\ref{cminimax}) if we set $f(x)=\delta_C(G(x))$ and $g(y)=\delta_D(H(y))$.
\end{example}
\quad \, The study of algorithms for solving convex-concave minimax problems of the form (\ref{cminimax}) is always active.
For the case when $K$ is a bilinear function, there are many publications about constructing and analyzing numerical algorithms for the minimax problem. The first work was due to  Arrow {\it et al.} \cite{Arrow58}, where they proposed an alternating coordinate method, leaving the convergence unsolved.
Nemirovski \cite{Nem2004} considered the following minimax problem
$$
\min_x \max_{y \in Y} g(x) + x^T Ay + h^Ty,
$$
where $Y$ is a compact convex set and $g$ is a ${\cal C}^{1,1}$ convex function. He proposed a mirror-prox algorithm which returns an approximate saddle point
within the complexity of ${\rm O}(1/\varepsilon)$.
 Nestrov \cite{Nes2007} developed a dual extrapolation algorithm for solving
variational inequalities which owns the complexity bound  ${\rm O}(1/\varepsilon)$  for Lipschitz continuous operators and applied the algorithm to  bilinear matrix games.   Chen {\it et al.} \cite{ChenLan2014}
  presented a novel accelerated primal-dual (APD) method for solving this class  of minimax problems, and showed that the APD method achieves the same optimal rate of
convergence as Nesterov's smoothing technique.
 Chambolle and Pock \cite{CPock2011} proposed a first-order primal-dual algorithm and established the convergence of the algorithm. Later,
Chambolle and Pock \cite{CPock2016a} provided the  ergodic convergence rate and Chambolle and Pock \cite{CPock2016b} explored the rate of convergence for accelerated primal-dual algorithms.

For smooth convex-concave minimax problems when $K$ is not bilinear, there are some works recently. Mokhtari {\it et al.} \cite{Ozd2019a} proposed  algorithms admitting a unified analysis as approximations of
the classical proximal point method for solving saddle point problems. Mokhtari {\it et al.}  \cite{Ozd2019b} proved that the optimistic gradient and extra-gradient methods achieve a convergence
rate of ${\rm O}(1/k)$ for smooth convex-concave saddle point problems.

 Recently, Lin {\it et al.} \cite{LinJordan2020} announced that they solved a longstanding open question pertaining to the design of near-optimal first-order algorithms for smooth and strongly-convex-strongly-concave minimax problems by presenting an algorithm with $\tilde{{\rm O}}(\sqrt{\kappa_x\kappa_y})$ gradient complexity, matching the lower bound
up to logarithmic factors.
  
  As far as we know, numerical methods for the nonsmooth convex-concave minimax problem (\ref{cminimax}) are rare. The only paper that came into our view was written by  Valkonen \cite{Valkonen2014}, in which  a natural
extension of the modified primal-dual hybrid gradient method, originally for
bilinear $K$, due to Chambolle and Pock, was studied.

Averaging the weights of neural nets is a prohibitive approach in particular because the zero-sum
game that is defined by training one deep net against another is not a convex-concave zero-sum
game. Thus it seems essential to identify training algorithms that make the last iterate of the training
be very close to the equilibrium, rather than only the average. However, the averaging technique is a popular tool for achieving good complexity, most of the mentioned  papers (including
\cite{Nem2004}, \cite{Nes2007}, \cite{CPock2011},\cite{CPock2016a},
\cite{ChenLan2014} and \cite{Ozd2019b}) adopted  that averaging technique. In this paper, we will not adopt the averaging technique in our algorithm.

The paper is organized as follows. In Section 2, we provide some technique results about  minimizing the sum of  a convex function and a quadratic proximal term, which play important roles in the convergence analysis of mspACM.
In Section 3, we propose mild assumptions and prove the global convergence property of mspACM for solving Problem (\ref{cminimax}). In Section 4, the linear rate of convergence of mspACM is demonstrated under the assumption that the solution mapping is locally metrically subregular. In Section 5, we report some preliminary numerical results for  linear regression problems with the saddle point formulation and  for solving separate linearly  constrained nonsmooth convex-concave minimax problems.  Some discussions are made in the last section.
\section{Preliminary}\label{Sec2}
\setcounter{equation}{0}
\quad \, In this section, we give several  results about properties of minimizing the sum of  a convex function and a quadratic proximal term. The results will be used in the next section.
\begin{lemma}\label{d-convex}
Let $f: {\cal Z} \rightarrow \Re$ be twice differentiable  and ${\cal Z}$ be a finite-dimensional real Hilbert space with a scalar product $\langle \cdot,\cdot \rangle$ and its induced norm $\|\cdot\|$. Assume that there exists a self-adjoint operator
$\Sigma_f: {\cal Z}\rightarrow {\cal Z}$ such that ${\rm D}^2 f(z) \succeq \Sigma_f$ for any $z \in {\cal Z}$ (which means that ${\rm D}^2 f(z) - \Sigma_f$ is a positive-definite self-adjoint operator). Define
\begin{equation}\label{eq:2.1}
q(z',z)=f(z)-\left [f(z')+{\rm D}f(z')(z-z')+\frac{1}{2}\|z-z'\|^2_{\Sigma_f}\right].
\end{equation}
Then for any $z' \in {\cal Z}$, $q(z',z)$ is a convex function of $z$.
\end{lemma}

For a smooth convex optimization problem, the proximal point method often involves minimizing the sum of a smooth function and a quadratic proximal term. An important inequality will be established in the following lemma.
\begin{lemma}\label{prox-lem}
Let the assumptions of Lemma \ref{d-convex} be satisfied and $T: {\cal Z}\rightarrow {\cal Z}$ be  a positively semi-definite self-adjoint operator such that
$$
T+\Sigma_f \succ 0.
$$
Then for any $z^c\in {\cal Z}$, the problem
$$
\min f(z)+\displaystyle \frac{1}{2}\|z-z^c\|^2_T
$$
has a unique solution, denoted by $z^+$. Moreover, for any $z\in {\cal Z}$,
\begin{equation}\label{eq:2.2}
f(z)+\displaystyle \frac{1}{2}\|z-z^c\|^2_T-\displaystyle \frac{1}{2}\|z-z^+\|^2_{T+\Sigma_f}\geq
f(z^+)+\displaystyle \frac{1}{2}\|z^+-z^c\|^2_T.
\end{equation}
\end{lemma}
{\bf Proof}. Define
$$
f_c(z)=f(z)+\displaystyle \frac{1}{2}\|z-z^c\|^2_T.
$$
Then
$$
f_c(z)=\left[f(z^+)+{\rm D}f(z^+)(z-z^+)+\frac{1}{2}\|z-z^+\|^2_{\Sigma_f}\right ]+\displaystyle \frac{1}{2}\|z-z^c\|^2_T+q(z^+,z).
$$
Let
$$
q_0(z)=q(z^+,z)+{\rm D}f(z^+)(z-z^+)+\displaystyle \frac{1}{2}\|z-z^c\|^2_T-\displaystyle \frac{1}{2}\|z-z^+\|^2_T.
$$
Then $q_0$ is convex and
$$
f_c(z)=f(z^+)+\displaystyle \frac{1}{2}\|z-z^+\|^2_T+\frac{1}{2}\|z-z^+\|^2_{\Sigma_f}+q_0(z).
$$
Since ${\rm D}f_c(z^+)=0$, we have ${\rm D}q_0(z^+)=0$ and $q_0$ arrives its minimum value at $z^+$.  Then for every $z \in {\cal Z}$,
$
q_0(z) \geq q_0(z^+)= \displaystyle \frac{1}{2}\|z^+-z^c\|^2_T
$, which is equivalent to
$$
f(z)+\displaystyle \frac{1}{2}\|z-z^c\|^2_T-\displaystyle \frac{1}{2}\|z-z^+\|^2_{T+\Sigma_f}\geq f(z^+)+\displaystyle \frac{1}{2}\|z^+-z^c\|^2_T.
$$
The proof is completed. \hfill $\Box$

Surprisingly, the inequality in Lemma \ref{prox-lem} is still valid when the convex function is the sum of a smooth convex function and a nonsmooth convex function. This result is given in the following lemma.
\begin{lemma}\label{prox-lem2}
Let the assumptions of Lemma \ref{d-convex} be satisfied and $T: {\cal Z}\rightarrow {\cal Z}$ be  a positively semi-definite self-adjoint operator such that
$$
T+\Sigma_f \succ 0.
$$
Then for any proper lower semicontinuous convex function $\psi:{\cal Z}\rightarrow \overline \Re$ and any $z^c\in {\cal Z}$, the problem
\begin{equation}\label{eq:composite}
\min \psi(z)+f(z)+\displaystyle \frac{1}{2}\|z-z^c\|^2_T
\end{equation}
has a unique solution, denoted by $z^+$. Moreover, for any $z\in {\cal Z}$,
\begin{equation}\label{eq:2.4}
\psi(z)+f(z)+\displaystyle \frac{1}{2}\|z-z^c\|^2_T-\displaystyle \frac{1}{2}\|z-z^+\|^2_{T+\Sigma_f}\geq
\psi(z^+)+f(z^+)+\displaystyle \frac{1}{2}\|z^+-z^c\|^2_T.
\end{equation}
\end{lemma}
{\bf Proof}. Define
$$
f_c(z)=\psi(z)+f(z)+\displaystyle \frac{1}{2}\|z-z^c\|^2_T.
$$
Then
$$
f_c(z)=\left[f(z^+)+{\rm D}f(z^+)(z-z^+)+\frac{1}{2}\|z-z^+\|^2_{\Sigma_f}\right]+\displaystyle \frac{1}{2}\|z-z^c\|^2_T+q(z^+,z)+\psi (z).
$$
Let
$$
q_0(z)=\psi (z)+q(z^+,z)+{\rm D}f(z^+)(z-z^+)+\displaystyle \frac{1}{2}\|z-z^c\|^2_T-\displaystyle \frac{1}{2}\|z-z^+\|^2_T.
$$
Then $q_0$ is convex and
$$
f_c(z)=f(z^+)+\displaystyle \frac{1}{2}\|z-z^+\|^2_T+\frac{1}{2}\|z-z^+\|^2_{\Sigma_f}+q_0(z).
$$
Since $0 \in \partial f_c(z^+)$, we have $0 \in \partial q_0(z^+)$ and $q_0$ arrives its minimum value at $z^+$.  Then for every $z \in {\cal Z}$,
$
q_0(z) \geq q_0(z^+)= \psi(z^+)+\displaystyle \frac{1}{2}\|z^+-z^c\|^2_T
$, which is equivalent to
$$
\psi(z)+f(z)+\displaystyle \frac{1}{2}\|z-z^c\|^2_T-\displaystyle \frac{1}{2}\|z-z^+\|^2_{T+\Sigma_f}\geq \psi(z^+)+f(z^+)+\displaystyle \frac{1}{2}\|z^+-z^c\|^2_T.
$$
The proof is completed. \hfill $\Box$

In the next section, we will find that the inequality (\ref{eq:2.4}) plays an important role in establishing the convergence of the proposed algorithm.
\section{The Algorithm and Global convergence}
\setcounter{equation}{0}

To begin with, we introduce some notation. Let ${\cal Z} ={\cal X}\times {\cal Y}$ and denote $z=(x,y)$, $z'=(x',y')$, $z^k=(x^k,y^k)$ and $z^{k+1/2}=(x^{k+1/2},y^{k+1/2})$. Let $K$ be a smooth convex-concave function on  an open set ${\cal O} \supset {\rm dom}\,f \times {\rm dom}\,g$; {\it i.e.}, for each $ (x,y)\in {\cal O}$, $K(\cdot,y)$ and $-K(x,\cdot)$ are smooth convex functions.
Let $\widehat \Sigma_f$ and $\widehat \Sigma_g$  be self-adjoint and positive semidefinite linear operators. We define, for any $(x,y),(x',y')\in {\cal O}$,
$$
\widehat K(x,y;x',y')=K(x',y')+\langle {\rm D}_x K(x',y'),x-x'\rangle+\langle {\rm D}_y K(x',y'),y-y'\rangle
+\displaystyle \frac{1}{2}\|x-x'\|^2_{\widehat \Sigma_f}-\displaystyle \frac{1}{2}\|y-y'\|^2_{\widehat \Sigma_g}
$$
and
$$
\widehat L(x,y;x',y')=f(x)+\widehat K(x,y;x',y')-g(y).
$$

We propose a majorized semi-proximal alternating coordinate algorithm (mspACM) for solving Problem (\ref{cminimax}) as below:
\begin{algorithm}\label{alg-1}
{\small	 {\rm (mspACM)}}
\begin{description}
\item[\bf Step 0.] Input $z^0=(x^0,y^0)\in {\rm dom}{\,f}\times{\rm dom}{\,g}$. Set $k:=0$.
\item[\bf Step 1.] Compute $z^{k+1/2}=(x^{k+1/2},y^{k+1/2})$ and $z^{k+1}=(x^{k+1},y^{k+1})$ by
$$
\left
\{
\begin{array}{ll}
 x^{k+1/2}&= \arg\min\limits_{x\in\mathcal{X}}~ \sigma \widehat L(x,y^k;z^k)+\displaystyle\frac{1}{2}\|x-x^k\|^2_{\mathcal{S}}\\
   &= \arg\min\limits_{x\in\mathcal{X}}~ \sigma \left[f(x)+\widehat K(x,y^k;z^k)-g(y^k)\right]+\displaystyle\frac{1}{2}\|x-x^k\|^2_{\mathcal{S}},\\
    y^{k+1/2}&= \arg\min\limits_{y \in\mathcal{Y}} -\sigma \widehat L(x^k,y;z^k)+\displaystyle\frac{1}{2}\|y-y^k\|^2_{\mathcal{T}}\\
   &= \arg\min\limits_{y\in\mathcal{Y}}~ \sigma \left[-f(x^k)-\widehat K(x^k,y;z^k)+g(y)\right]+\displaystyle\frac{1}{2}\|y-y^k\|^2_{\mathcal{T}},\\
x^{k+1} &= \arg\min\limits_{x\in \mathcal{X}} \sigma \widehat L(x,y^{k+1/2};z^{k+1/2})+\displaystyle\frac{1}{2}\|x-x^k\|^2_{\mathcal{S}}\\
   &= \arg\min\limits_{x\in\mathcal{X}} \sigma \left[f(x)+\widehat K(x,y^{k+1/2};z^{k+1/2})-g(y^{k+1/2})\right]+\displaystyle\frac{1}{2}\|x-x^k\|^2_{\mathcal{S}},\\
    y^{k+1}&= \arg\min\limits_{y\in\mathcal{Y}}~ -\sigma \widehat L(x^{k+1/2},y;z^{k+1/2})+\displaystyle\frac{1}{2}\|y-y^k\|^2_{\mathcal{T}}\\
   &= \arg\min\limits_{y\in\mathcal{Y}}~ \sigma \left[-f(x^{k+1/2})-\widehat K(x^{k+1/2},y;z^{k+1/2})+g(y)\right]+\displaystyle\frac{1}{2}\|y-y^k\|^2_{\mathcal{T}}.
\end{array}
\right.
$$
\item[\bf Step 2.] If a termination criterion is not met, set $k:=k+1$ and go to Step 1.
\end{description}
\end{algorithm}

The motivation for proposing the above algorithm comes from the following observations:
\begin{itemize}
\item[{\rm (i)}] The reason for using $\widehat K(\cdot,\cdot;z)$ in stead of $K(\cdot,\cdot)$ is that this makes
the subproblems for determining $z^{k+1/2}$ and $z^{k+1}$ easily be solved, especially when $K$ is a complicated smooth convex function. Furthermore, the subproblems for determining $z^{k+1/2}$ and $z^{k+1}$ may have explicit solutions when $f$ and $g$ are simple convex functions.
\item[{\rm (ii)}] The use of the semi-proximal terms ({\it i.e.}, ${\cal S}$ and ${\cal T}$ are only required to be positively semidefinite)  leaves the user a freedom to choose ${\cal S}$ and ${\cal T}$ so that the subproblems for determining $z^{k+1/2}$ and $z^{k+1}$ are well-conditioned or are easily solved.
\end{itemize}

For analyzing the global convergence of Algorithm \ref{alg-1}, we need the following mild assumptions about functions  in Problem (\ref{cminimax}).
\begin{assumption}\label{K-cc}
Let $K$ be a smooth convex-concave function on  an open set ${\cal O} \supset {\rm dom}\,f \times {\rm dom}\,g$; {\it i.e.,} for each $ (x,y)\in {\cal O}$, $K(\cdot,y)$ and $-K(x,\cdot)$ are smooth convex functions.
Suppose there exist self-adjoint and positive semidefinite linear operators $\widehat \Sigma_f$ and $\widehat \Sigma_g$ such that for any $(x,y),(x',y')\in {\cal O}$,
\begin{eqnarray}
 K(x,y) & \leq & K(x',y)+\langle {\rm D}_x K(x',y),x-x'\rangle+\frac{1}{2}\|x-x'\|^2_{\widehat \Sigma_f},\label{def-LB-f}\\
  -K(x,y) & \leq & -K(x,y')-\langle {\rm D}_y K(x,y'),y-y'\rangle+\frac{1}{2}\|y-y'\|^2_{\widehat \Sigma_g}\label{def-LB-g}.
\end{eqnarray}
\end{assumption}
For convenience, we introduce a linear operator $\widehat \Sigma:{\cal Z}\rightarrow {\cal Z}$ by
\begin{equation}\label{hatSigma}
\widehat \Sigma (x,y)=(\widehat \Sigma_f x,\widehat \Sigma_g y).
\end{equation}
\begin{assumption}\label{assu-K}
Suppose that $K$ is continuously differentiable on an open set ${\cal O} \supset {\rm dom}\,f \times {\rm dom}\,g$  and the derivative mapping ${\rm D} K$ is Lipschitz continuous with constant $\eta_0>0$; {\it i.e.},
$$
\|{\rm D}K(x,y)-{\rm D}K(x',y')\|\leq \eta_0\|(x,y)-(x',y')\|, \ \  (x',y),(x,y) \in {\cal O}.
$$
\end{assumption}
\begin{assumption}\label{ass-o}
The set of saddle points of $L$ over ${\rm dom}\,f \times {\rm dom}\,g$ is nonempty; {\i.e.},
 $\overline \Omega\ne \emptyset$, where  $\overline \Omega$ is defined by
 $$
\overline \Omega=\left\{(x,y) \in {\cal Z}: L(x,y')\leq L(x,y) \leq L(x',y),\forall (x',y')\in {\rm dom}\,f \times {\rm dom}\,g\right \}.
 $$
 \end{assumption}
 Define for $z=(x,y)$, $z'=(x',y')$,
$$
\Phi_K(z',z)=K(x',y)-K(x,y'),\quad  \Phi_L(z',z)=L(x',y)-L(x,y').
$$
Then we have for $z',z\in  {\rm dom}\,f \times {\rm dom}\,g$ that
$$
\Phi_K(z,z)=0,\ \Phi_L(z,z)=0,\ \Phi_K(z',z)+\Phi_K(z,z')=0,\ \Phi_L(z',z)+\Phi_L(z,z')=0.
$$
\begin{proposition}\label{4diff}
Let Assumption \ref{assu-K} be satisfied. Then for $z',z\in  {\rm dom}\,f \times {\rm dom}\,g$ and $h',h \in {\cal X} \times {\cal Y}$ such that $z'+h',z+h\in {\rm dom}\,f \times {\rm dom}\,g$, one has that
$$
\|[\Phi_K(z'+h',z+h)-\Phi_K(z',z+h)]-[\Phi_K(z'+h',z)-\Phi_K(z',z)]\|\leq \eta_0 \|h'\|\|h\|
$$
 and
  $$
  \begin{array}{l}
 \|[\Phi_L(z'+h',z+h)-\Phi_L(z',z+h)]-[\Phi_L(z'+h',z)-\Phi_L(z',z)]\\[4pt]
 \quad =\|[\Phi_K(z'+h',z+h)-\Phi_K(z',z+h)]-[\Phi_K(z'+h',z)-\Phi_K(z',z)]\|
 \end{array}
 $$
 so that
 $$
\|[\Phi_L(z'+h',z+h)-\Phi_L(z',z+h)]-[\Phi_L(z'+h',z)-\Phi_L(z',z)]\|\leq \eta_0 \|h'\|\|h\|.
$$
\end{proposition}
Furthermore, define for $z=(x,y)$, $z'=(x',y')$ and $z''=(x'',y'')$,
$$
\widehat \Phi_K(z',z'';z)=\widehat K(x',y'';z)-\widehat K(x'',y';z),\quad  \widehat \Phi_L(z',z'';z)=\widehat L(x',y'';z)-\widehat L(x'',y';z).
$$
Then we have for $z',z'',z\in  {\rm dom}\,f \times {\rm dom}\,g$ that
$$
\widehat \Phi_K(z',z';z)=0,\Phi_L(z',z';z)=0,\widehat \Phi_K(z',z'';z)+\widehat \Phi_K(z'',z';z)=0,\widehat \Phi_L(z',z'';z)+\widehat \Phi_L(z'',z';z)=0.
$$
\begin{proposition}\label{4diffa}
Let Assumption \ref{assu-K} be satisfied. Then for $z\in  {\rm dom}\,f \times {\rm dom}\,g$ and $h',h \in {\cal X} \times {\cal Y}$ such that $z+h+h',z+h\in {\rm dom}\,f \times {\rm dom}\,g$, one has that
$$
\|[\widehat\Phi_K(z+h+h',z;z)-\widehat\Phi_K(z+h,z;z)]-[\widehat\Phi_K(z+h+h',z+h;z+h)-\widehat\Phi_K(z+h,z+h;z+h)]\|\leq \widehat\eta_0 \|h\|\|h'\|
$$
 and
  $$
  \begin{array}{l}
 \|[\widehat\Phi_L(z+h+h',z;z)-\widehat\Phi_L(z+h,z;z)]-[\widehat\Phi_L(z+h+h',z+h;z+h)-\widehat\Phi_L(z+h,z+h;z+h)]\|\\[4pt]
 \quad =\|[\widehat\Phi_K(z+h+h',z;z)-\widehat\Phi_K(z+h,z;z)]-[\widehat\Phi_K(z+h+h',z+h;z+h)-\widehat\Phi_K(z+h,z+h;z+h)]\|
 \end{array}
 $$
 so that
 $$
\|[\widehat\Phi_L(z+h+h',z;z)-\widehat\Phi_L(z+h,z;z)]-[\widehat\Phi_L(z+h+h',z+h;z+h)-\widehat\Phi_L(z+h,z+h;z+h)]\|
\leq\widehat \eta_0 \|h\|\|h'\|,
$$
where
$$
\widehat\eta_0=\|\widehat \Sigma\|+\eta_0.
$$
\end{proposition}
{\bf Proof}. Let $h=(h_x,h_y)\in {\cal Z}$. Define an operation $\tilde{\rm D}$ by
$$
\widetilde{\rm D}K(z)h=\left({\rm D}_xK(z)h_x, -{\rm D}_yK(z)h_y\right).$$
Since $\widehat\Phi_K(z+h,z+h;z+h)=0$, we only need to consider the other three terms. For $h=(h_x,h_y)$ and $h'=(h'_x,h'_y),$
\begin{eqnarray}
\widehat\Phi_K(z+h+h',z;z)&=&\widehat K(x+h_x+h'_x,y;z)-\widehat K(x,y+h_y+h'_y;z)\nonumber \\
&=&{\rm D}_x K(z)(h_x+h'_x)-{\rm D}_y K(z)(h_y+h'_y)+\displaystyle \frac{1}{2}\|h_x+h'_x\|_{\widehat \Sigma_f}^2+
\displaystyle \frac{1}{2}\|h_y+h'_y\|_{\widehat \Sigma_g}^2 \nonumber \\
&=&\widetilde{\rm D} K(z)(h+h')+\displaystyle \frac{1}{2}\|h+h'\|^2_{\widehat \Sigma}.\label{f1}
\end{eqnarray}
Similarly, we can get that
\begin{equation}\label{f2}
\begin{array}{ll}
\widehat\Phi_K(z+h,z;z)=\widetilde{\rm D} K(z)(h)+\displaystyle \frac{1}{2}\|h\|^2_{\widehat \Sigma}.
\end{array}
\end{equation}
and
\begin{equation}\label{f3}
\begin{array}{ll}
\widehat\Phi_K(z+h+h',z+h;z+h)
=\widetilde{\rm D} K(z+h)(h')+\displaystyle \frac{1}{2}\|h'\|^2_{\widehat \Sigma}.
\end{array}
\end{equation}
Combing (\ref{f1}),(\ref{f2}) and (\ref{f3}), we obtain
$$
\begin{array}{l}
|[\widehat\Phi_K(z+h+h',z;z)-\widehat\Phi_K(z+h,z;z)]-[\widehat\Phi_K(z+h+h',z+h;z+h)-\widehat\Phi_K(z+h,z+h;z+h)]|\\[6pt]
=|[\widehat\Phi_K(z+h+h',z;z)-\widehat\Phi_K(z+h,z;z)]-[\widehat\Phi_K(z+h+h',z+h;z+h)|\\[6pt]
=|[\widetilde{\rm D} K(z)-\widetilde{\rm D} K(z+h)](h')+\displaystyle \frac{1}{2}\|h+h'\|^2_{\widehat \Sigma}-
\displaystyle \frac{1}{2}\|h\|^2_{\widehat \Sigma}-\displaystyle \frac{1}{2}\|h'\|^2_{\widehat \Sigma}|\\[4pt]
=|[\widetilde{\rm D} K(z)-\widetilde{\rm D} K(z+h)](h')+\langle\widehat \Sigma h, h' \rangle|\\[4pt]
\leq \|[\widetilde{\rm D} K(z)-\widetilde{\rm D} K(z+h)]\|\|h'\|+\|\widehat \Sigma \|\|h\|\|h'\|\\[4pt]
\leq (\eta_0+\|\widehat \Sigma \|)\|h\|\|h'\|=\widehat\eta_0\|h\|\|h'\|.
\end{array}
$$
The proof is completed.\hfill $\Box$\\

Let $\sigma>0$  be a given parameter. Let $\mathcal{S}$ and $\mathcal{T}$ be given self-adjoint linear operators satisfying
\begin{equation}\label{cond-2}
\begin{array}{c}
\mathcal{S}\succeq 0,\ \ \mathcal{T} \succeq 0,\ \ \widehat \Sigma_f+\mathcal{S}\succ 0,\ \ \widehat \Sigma_g+\mathcal{T}\succ 0.
\end{array}
\end{equation}
Define a linear operator $\Theta:{\cal Z}\rightarrow {\cal Z}$ by
\begin{equation}\label{Theta}
\Theta (x,y)=({\cal S}x, {\cal T}y).
\end{equation}
Then the formulas for $(x^{k+1/2},y^{k+1/2})$ and $(x^{k+1},y^{k+1})$ in Step 1 of Algorithm \ref{alg-1} can be written as
\begin{equation}\label{key-o}
\begin{array}{ll}
z^{k+1/2}&=\arg\min\limits_{z\in\mathcal{Z}}\left\{\displaystyle \frac{1}{2}\|z-z^k\|^2_{\Theta}+\sigma \widehat\Phi_L(z,z^k;z^k)\right\},\\[4pt]
z^{k+1}&=\arg\min\limits_{z\in\mathcal{Z}}\left\{\displaystyle \frac{1}{2}\|z-z^k\|^2_{\Theta}+\sigma \widehat\Phi_L(z,z^{k+1/2};z^{k+1/2})\right\}.
\end{array}
\end{equation}
 Noting that $\Phi_L(z,z^{k+1/2};z^{k+1/2})$ and $\Phi_L(z,z^{k};z^k)$ are convex
and problems in (\ref{key-o}) have the same  structure as problem
(\ref{eq:composite}), we may use Lemma \ref{prox-lem2} to estimate $\|z-z^{k+1}\|_{\widehat\Sigma+\Theta}$ and $\|z-z^{k+1/2}\|_{\widehat\Sigma+\Theta}$. Specifically, by Lemma \ref{prox-lem2}, we know from the definition of $z^{k+1}$ that for any $z \in {\rm dom}{\,f}\times{\rm dom}{\,g}$,
\begin{equation}\label{a1}
\displaystyle \frac{1}{2}\|z-z^k\|^2_{\Theta} +\sigma \widehat \Phi_L(z,z^{k+1/2};z^{k+1/2})-\displaystyle \frac{1}{2}\|z-z^{k+1}\|^2_{\widehat\Sigma+\Theta}
\geq \displaystyle \frac{1}{2}\|z^{k+1}-z^k\|^2_{\Theta} +\sigma \widehat\Phi_L(z^{k+1},z^{k+1/2};z^{k+1/2}).
\end{equation}
Similarly, we have from the definition of $z^{k+1/2}$ that for any $z \in {\rm dom}{\,f}\times{\rm dom}{\,g}$,
\begin{equation}\label{a2}
\displaystyle \frac{1}{2}\|z-z^k\|^2_{\Theta} +\sigma \widehat\Phi_L(z,z^{k};z^k)-\displaystyle \frac{1}{2}\|z-z^{k+1/2}\|^2_{\widehat\Sigma+\Theta}
\geq \displaystyle \frac{1}{2}\|z^{k+1/2}-z^k\|^2_{\Theta} +\sigma \widehat\Phi_L(z^{k+1/2},z^{k};z^k).
\end{equation}

We first establish the relation between $\|z^{k+1}-z^{k+1/2}\|_{\Theta+\widehat\Sigma}$ and $\|z^{k+1/2}-z^{k}\|_{\Theta+\widehat\Sigma}$ in the following proposition.
\begin{proposition}\label{relationB}
Let Assumption \ref{K-cc} and Assumption \ref{assu-K} be satisfied, ${\cal S}$ and ${\cal T}$ satisfy (\ref{cond-2}) or equivalently $\Theta \succeq 0$ and $\Theta+\widehat \Sigma\succ 0$.  Let $\{z^k=(x^k,y^k)\},\{z^{k+1/2}=(x^{k+1/2},y^{k+1/2})\}$  be generated by Algorithm
\ref{alg-1}. Then
\begin{equation}\label{eq:relation-z}
\|z^{k+1}-z^{k+1/2}\|_{\Theta+\widehat \Sigma} \leq \vartheta(\sigma) \|z^{k+1/2}-z^{k}\|_{\Theta+\widehat \Sigma},
\end{equation}
where
$$
\vartheta(\sigma)=\displaystyle \frac{\sigma \widehat\eta_0}{\lambda_{\min}(\Theta+\widehat \Sigma)}.
$$
\end{proposition}
{\bf Proof}.
Setting $z=z^{k+1/2}$ in (\ref{a1}), we obtain
\begin{equation}\label{a1a}
\displaystyle \frac{1}{2}\|z^{k+1/2}-z^k\|^2_{\Theta} +\sigma \widehat\Phi_L(z^{k+1/2},z^{k+1/2};z^{k+1/2})-\displaystyle \frac{1}{2}\|z^{k+1/2}-z^{k+1}\|^2_{\widehat \Sigma+\Theta}
\geq \displaystyle \frac{1}{2}\|z^{k+1}-z^k\|^2_{\Theta} +\sigma \widehat\Phi_L(z^{k+1},z^{k+1/2};z^{k+1/2}).
\end{equation}
Setting $z=z^{k+1}$ in (\ref{a2}), we obtain
\begin{equation}\label{a2a}
\displaystyle \frac{1}{2}\|z^{k+1}-z^k\|^2_{\Theta} +\sigma \widehat\Phi_L(z^{k+1},z^{k};z^k)-\displaystyle \frac{1}{2}\|z^{k+1}-z^{k+1/2}\|^2_{\widehat \Sigma+\Theta}
\geq \displaystyle \frac{1}{2}\|z^{k+1/2}-z^k\|^2_{\Theta} +\sigma \widehat\Phi_L(z^{k+1/2},z^{k};z^k).
\end{equation}
Summing (\ref{a1a}) and (\ref{a2a}), we get from Proposition \ref{4diffa} (with $z=z^k$, $h=z^{k+1/2}-z^k$, $h'=z^{k+1}-z^{k+1/2}$)
 that
\begin{eqnarray}
&& \!\!\!\!\! \!\!\!\!\!\!\!\|z^{k+1}-z^{k+1/2}\|^2_{\widehat \Sigma+\Theta} \nonumber\\
 &\leq& \!\!\!\!\sigma\{\widehat\Phi_L(z^{k+1/2},z^{k+1/2};z^{k+1/2})-\widehat\Phi_L(z^{k+1},z^{k+1/2};z^{k+1/2})-
[\widehat\Phi_L(z^{k+1/2},z^{k};z^k)-\widehat\Phi_L(z^{k+1},z^{k};z^k)]\}\nonumber\\
&=&\!\!\!\!\sigma\{\widehat\Phi_K(z^{k+1/2},z^{k+1/2};z^{k+1/2})-\widehat\Phi_K(z^{k+1},z^{k+1/2};z^{k+1/2})-
[\widehat\Phi_K(z^{k+1/2},z^{k};z^k)-\widehat\Phi_K(z^{k+1},z^{k};z^k)]\}\nonumber\\
&\leq& \!\!\!\!\sigma \widehat\eta_0 \|z^{k+1}-z^{k+1/2}\|\|z^{k+1/2}-z^k\|. \label{a3}
\end{eqnarray}
Therefore, in view of
$$
\sqrt{\lambda_{\min}(\widehat \Sigma+\Theta)}\|z\| \leq \|z\|_{\widehat \Sigma+\Theta}, \forall z \in {\cal Z},
$$
we obtain
 the desired result from the last expression of (\ref{a3}). The proof is completed. \hfill $\Box$

 Now we provide the main result about the global convergence of Algorithm \ref{alg-1}.
 \begin{theorem}\label{gconv}
Let Assumption \ref{K-cc}, Assumption \ref{assu-K} and Assumption \ref{ass-o} be satisfied and ${\cal S}$ and ${\cal T}$ satisfy (\ref{cond-2}) or equivalently $\Theta \succeq 0$ and $\Theta+\widehat \Sigma\succ 0$.  Consider the sequences $\{z^k=(x^k,y^k)\}$ and $\{z^{k+1/2}=(x^{k+1/2},y^{k+1/2})\}$ generated by Algorithm \ref{alg-1}.
Suppose that $\sigma$, ${\cal S}$ and ${\cal T}$ satisfy
\begin{equation}\label{eq:par}
0< \sigma <\displaystyle \min\left\{ \displaystyle \frac{\lambda_{\min}(\widehat \Sigma+\Theta)}{\sqrt{2}\widehat\eta_0},\displaystyle \frac{1}{2}\right\},\ \ \Theta \succ (\widehat \eta_0+2)\sigma I.
\end{equation}
If $\widehat \Sigma$ satisfies
\begin{equation}\label{eq:hatsigma}
\widehat \Sigma \succ \displaystyle \frac{6\eta_0\sigma}{1-2\sigma} I,
\end{equation}
then $\{z^k=(x^k,y^k)\}$ converges monotonically with respect to some
norm to an element of $\overline \Omega$.
\end{theorem}
{\bf Proof}. Choose an element $z^* \in \overline \Omega$. Setting $z=z^*$ in (\ref{a1}), we obtain
\begin{equation}\label{a1aa}
\displaystyle \frac{1}{2}\|z^*-z^k\|^2_{\Theta} +\sigma \widehat \Phi_L(z^*,z^{k+1/2};z^{k+1/2})-\|z^*-z^{k+1}\|^2_{\widehat \Sigma+\Theta}
\geq \displaystyle \frac{1}{2}\|z^{k+1}-z^k\|^2_{\Theta} +\sigma \widehat\Phi_L(z^{k+1},z^{k+1/2};z^{k+1/2}).
\end{equation}
Setting $z=z^{k+1}$ in (\ref{a2}), we obtain
\begin{equation}\label{a2aa}
\displaystyle \frac{1}{2}\|z^{k+1}-z^k\|^2_{\Theta} +\sigma \widehat \Phi_L(z^{k+1},z^{k};z^k)-\|z^{k+1}-z^{k+1/2}\|^2_{\widehat \Sigma+\Theta}
\geq \displaystyle \frac{1}{2}\|z^{k+1/2}-z^k\|^2_{\Theta} +\sigma \widehat\Phi_L(z^{k+1/2},z^{k};z^k).
\end{equation}
Summing (\ref{a1aa}) and (\ref{a2aa}), we have
\begin{eqnarray}
&& \!\!\!\!\!\!\!\!\!\!\|z^k-z^*\|_{\Theta}^2\nonumber \\
& \geq & \!\!\!\! \|z^{k+1}-z^*\|_{\widehat \Sigma+\Theta}^2+\|z^{k+1}-z^{k+1/2}\|^2_{\widehat \Sigma+\Theta}+\|z^{k+1/2}-z^k\|^2_{\Theta}\nonumber \\
&&\!\!\!\!+
2\sigma[\widehat \Phi_L(z^{k+1},z^{k+1/2};z^{k+1/2})-\widehat\Phi_L(z^*,z^{k+1/2};z^{k+1/2})+\widehat\Phi_L(z^{k+1/2},z^{k};z^k)-\widehat\Phi_L(z^{k+1},z^{k};z^k)]\nonumber \\
&=&\!\!\!\!\|z^{k+1}-z^*\|_{\widehat \Sigma+\Theta}^2+\|z^{k+1}-z^{k+1/2}\|^2_{\widehat \Sigma+\Theta}+\|z^{k+1/2}-z^k\|_{\Theta}^2\nonumber \\
& &\!\!\!\!+
2\sigma\left[\widehat\Phi_L(z^{k+1},z^{k+1/2};z^{k+1/2})-\widehat\Phi_L(z^{k+1/2},z^{k+1/2};z^{k+1/2})+\widehat\Phi_L(z^{k+1/2},z^{k};z^k)
-\widehat\Phi_L(z^{k+1},z^{k};z^k)\right]\nonumber \\
& &\!\!\!\!+2\sigma\left[\widehat\Phi_L(z^{k+1/2},z^{k+1/2};z^{k+1/2})-\widehat\Phi_L(z^*,z^{k+1/2};z^{k+1/2})\right]. \label{a5}
\end{eqnarray}
Noting that
$$
\widehat \Phi_L(z^{k+1/2},z^{k+1/2};z^{k+1/2})=0,
$$
we need to estimate the term $-\widehat \Phi_L(z^*,z^{k+1/2};z^{k+1/2})$. In fact, we have that
\begin{eqnarray}
&& \!\!\!\!\!\!\!\!\!\! -\widehat \Phi_L(z^*,z^{k+1/2};z^{k+1/2}) \nonumber \\
&=&\!\!\!\!\widehat L(x^{k+1/2},y^*;z^{k+1/2})-\widehat L(x^*,y^{k+1/2};z^{k+1/2})\nonumber\\
&=&\!\!\!\! f(x^{k+1/2})+\widehat K(x^{k+1/2},y^*;z^{k+1/2})-g(y^*)-\left[f(x^*)+\widehat K(x^*,y^{k+1/2};z^{k+1/2})-g(y^{k+1/2})\right]\nonumber\\
&=&\!\!\!\! f(x^{k+1/2})-g(y^*)+K(z^{k+1/2})+\langle{\rm D}_yK(z^{k+1/2}),y^*-y^{k+1/2}\rangle-\displaystyle \frac{1}{2}\|y^*-y^{k+1/2}\|_{\widehat \Sigma_g}^2\nonumber\\
&&\!\!\!\! -\left[f(x^*)-g(y^{k+1/2})+K(z^{k+1})+\langle{\rm D}_xK(z^{k+1/2}),x^*-x^{k+1/2}\rangle\right]
  -\displaystyle \frac{1}{2}\|x^*-x^{k+1/2}\|_{\widehat \Sigma_f}^2 \nonumber\\
&=&\!\!\!\! f(x^{k+1/2})-g(y^*)-f(x^*)+g(y^{k+1/2})  +\langle{\rm D}_yK(z^{k+1/2}),y^*-y^{k+1/2}\rangle \nonumber \\
& &\!\!\!\! +\langle{\rm D}_xK(z^{k+1/2}),x^{k+1/2}-x^*\rangle -\displaystyle \frac{1}{2}\|y^*-y^{k+1/2}\|_{\widehat \Sigma_g}^2
 -\displaystyle \frac{1}{2}\|x^*-x^{k+1/2}\|_{\widehat \Sigma_f}^2 \nonumber \\
&=&\!\!\!\! f(x^{k+1/2})-g(y^*)-f(x^*)+g(y^{k+1/2}) +\langle{\rm D}_xK(z^{k+1/2}),x^{k+1/2}-x^*\rangle \nonumber \\
&&\!\!\!\! +\langle{\rm D}_yK(z^{k+1/2}),y^*-y^{k+1/2}\rangle-\displaystyle \frac{1}{2}\|z^{k+1/2}-z^*\|_{\widehat \Sigma}^2.\label{eq:b1}
\end{eqnarray}

Since $(x^*,y^*)\in \overline \Omega$ and $(x^{k+1/2},y^{k+1/2})\in {\rm dom}\,f \times {\rm dom}\,g$, we have
\begin{equation}\label{eq:b2}
\begin{array}{l}
f(x^{k+1/2})-g(y^*)-f(x^*)+g(y^{k+1/2})\\[4pt]
=\left[f(x^{k+1/2})-g(y^*)+K(x^{k+1/2},y^*)\right]-L(x^*,y^*)\\[4pt]
\quad +L(x^*,y^*)-\left[f(x^*)+g(y^{k+1/2})+K(x^*,y^{k+1/2})\right]
   +K(x^*,y^{k+1/2})-K(x^{k+1/2},y^*)\\[4pt]
=\left[L(x^{k+1/2},y^*)-L(x^*,y^*)\right]+\left[L(x^*,y^*)-L(x^*,y^{k+1/2})\right]+K(x^*,y^{k+1/2})-K(x^{k+1/2},y^*)\\[4pt]
\geq K(x^*,y^{k+1/2})-K(x^*,y^*)-\left[K(x^{k+1/2},y^*)-K(x^*,y^*)\right]\\[4pt]
=\left\langle \int_0^1{\rm D}_y K(x^*,y^k(t)){\rm d}t, y^{k+1/2}-y^*\right\rangle-\left\langle \int_0^1{\rm D}_x K(x^k(t),y^*){\rm d}t, x^{k+1/2}-x^*\right\rangle,
\end{array}
\end{equation}
where $y^k(t)=y^*+t(y^{k+1/2}-y^*)$ and $x^k(t)=x^*+t(x^{k+1/2}-x^*)$ for $t \in [0,1]$.
It follows from (\ref{eq:b2}) that
\begin{eqnarray}
&&\!\!\!\!\!\!\!\!\!\! f(x^{k+1/2})-g(y^*)-f(x^*)+g(y^{k+1/2})
+\langle{\rm D}_xK(z^{k+1/2}),x^{k+1/2}-x^*\rangle
 +\langle{\rm D}_yK(z^{k+1/2}),y^*-y^{k+1/2}\rangle \nonumber\\
&\geq&\!\!\!\! \left\langle\displaystyle \int_0^1[{\rm D}_y K(x^*,y^k(t))-{\rm D}_yK(z^k)]{\rm d}t, y^{k+1/2}-y^*\right\rangle
 -\left\langle \displaystyle\int_0^1[{\rm D}_x K(x^k(t),y^*)-{\rm D}_xK(z^k)]{\rm d}t, x^{k+1/2}-x^*\right\rangle \nonumber\\
& \geq&\!\!\!\! -\eta_0\displaystyle\int_0^1 \|(x^*,y^k(t))-z^{k+1/2}\|\|y^{k+1/2}-y^*\|{\rm d}t -\eta_0\displaystyle\int_0^1 \|(x^k(t),y^*)-z^{k+1/2}\|\|x^{k+1/2}-x^*\|{\rm d}t\nonumber\\
&\geq&\!\!\!\! -\eta_0\left\{\left[\|x^{k+1/2}-x^*\|+\|y^{k+1/2}-y^*\|/2\right]\|y^{k+1/2}-y^*\|\right. \nonumber\\
& &\!\!\!\! + \left.\left[ \|y^{k+1/2}-y^*\|+\|x^{k+1/2}-x^*\|/2\right]\|x^{k+1/2}-x^*\|\right\}\nonumber\\
&=&\!\!\!\!-\eta_0 \displaystyle\left( \frac{1}{2}\|z^{k+1/2}-z^*\|^2+2 \|x^{k+1/2}-x^*\|\|y^{k+1/2}-y^*\|\right). \label{eq:b3}
     \end{eqnarray}
Thus we can get from (\ref{a5}), (\ref{eq:b1}) and (\ref{eq:b3}) that
\begin{eqnarray}
\|z^k-z^*\|_{\Theta}^2 &\!\!\!\! \geq  &\!\!\!\! \|z^{k+1}-z^*\|_{\widehat \Sigma+\Theta}^2+\|z^{k+1}-z^{k+1/2}\|^2_{\widehat \Sigma+\Theta}+\|z^{k+1/2}-z^k\|_{\Theta}^2+
2\sigma\left[\widehat\Phi_L(z^{k+1},z^{k+1/2};z^{k+1/2}) \right.\nonumber \\[4pt]
&&\!\!\!\!-\left.\widehat\Phi_L(z^{k+1/2},z^{k+1/2};z^{k+1/2})+
\widehat\Phi_L(z^{k+1/2},z^{k};z^{k+1/2})-\widehat \Phi_L(z^{k+1},z^{k};z^k)\right] \nonumber\\[4pt]
&& \!\!\!\! -\sigma\|z^{k+1/2}-z^*\|_{\widehat \Sigma}^2-\eta_0\sigma \displaystyle\left( \|z^{k+1/2}-z^*\|^2+4 \|x^{k+1/2}-x^*\|\|y^{k+1/2}-y^*\|\right). \label{a6}
\end{eqnarray}
In view of Proposition \ref{4diffa},  we obtain from  (\ref{a6}) that
\begin{eqnarray}
&&\!\!\!\!\!\!\!\!\!\!\|z^k-z^*\|_{\Theta}^2\nonumber \\
 &\!\!\!\!\geq &\!\!\!\! \|z^{k+1}-z^*\|_{\widehat \Sigma+\Theta}^2+\|z^{k+1}-z^{k+1/2}\|^2_{\widehat \Sigma+\Theta}+\|z^{k+1/2}-z^k\|_{\Theta}^2-2\sigma\widehat \eta_0\|z^{k+1}-z^{k+1/2} \|\|z^{k+1/2}-z^{k}\|\nonumber\\[8pt]
&&\!\!\!\! -\sigma\|z^{k+1/2}-z^*\|_{\widehat \Sigma}^2-\eta_0\sigma \displaystyle\left( \|z^{k+1/2}-z^*\|^2+4 \|x^{k+1/2}-x^*\|\|y^{k+1/2}-y^*\|\right)\nonumber\\[8pt]
&\!\!\!\! \geq &\!\!\!\! \|z^{k+1}-z^*\|_{\widehat \Sigma+\Theta}^2+\|z^{k+1}-z^{k+1/2}\|^2_{\widehat \Sigma+\Theta}+\|z^{k+1/2}-z^k\|_{\Theta}^2-2\sigma\widehat \eta_0\|z^{k+1}-z^{k+1/2} \|\|z^{k+1/2}-z^{k}\|\nonumber\\[8pt]
&&\!\!\!\! -\sigma\|z^{k+1/2}-z^*\|_{\widehat \Sigma}^2-3\eta_0\sigma  \|z^{k+1/2}-z^*\|^2\nonumber\\[8pt]
&\!\!\!\!\geq &\!\!\!\! \|z^{k+1}-z^*\|_{\widehat \Sigma+\Theta}^2+\|z^{k+1}-z^{k+1/2}\|^2_{\widehat \Sigma+\Theta}+\|z^{k+1/2}-z^k\|_{\Theta}^2-\sigma\widehat \eta_0(\|z^{k+1}-z^{k+1/2} \|^2+\|z^{k+1/2}-z^{k}\|^2)\nonumber\\[8pt]
&&\!\!\!\! -2\sigma(\|z^{k+1/2}-z^k\|_{\widehat \Sigma}^2+\|z^{k}-z^*\|_{\widehat \Sigma}^2)
 -6\eta_0\sigma(\|z^{k+1}-z^{k+1/2}\|^2 + \|z^{k+1}-z^*\|^2). \label{a67}
\end{eqnarray}
Define
$$
\begin{array}{l}
G(\sigma)=\widehat \Sigma+\Theta-6\eta_0\sigma I,\\[4pt]
N(\sigma)=\widehat \Sigma+\Theta-(\widehat \eta_0+6\eta_0)\sigma I,\\[4pt]
H(\sigma)=\Theta-(\widehat \eta_0+2)\sigma I.
\end{array}
$$
The relation (\ref{a67}) implies that
\begin{equation}\label{ac7}
\|z^k-z^*\|_{\Theta+2\sigma\widehat \Sigma}^2  \geq  \|z^{k+1}-z^*\|_{G(\sigma)}^2+\|z^{k+1}-z^{k+1/2}\|^2_{N(\sigma)}\\[6pt]
+\|z^{k+1/2}-z^k\|_{H(\sigma)}^2.
\end{equation}
From (\ref{eq:par}), we know that $H(\sigma)$ is positively definite. From (\ref{eq:hatsigma}), we have that
$$
G(\sigma)=\widehat \Sigma+\Theta -6\eta_0 \sigma I
\succeq \Theta+2\sigma\widehat \Sigma
$$
and
$$
N(\sigma)=\widehat \Sigma+\Theta-(\widehat \eta_0+6\eta_0)\sigma I \succ
(\widehat \eta_0+2)\sigma I+\displaystyle \frac{6\eta_0\sigma}{1-2\sigma} I-(\widehat \eta_0+6\eta_0)\sigma I \succeq 2\sigma I.
$$
Thus both $G(\sigma)$ and $N(\sigma)$ are positively definite. Hence by (\ref{ac7}), we get that
\begin{equation}\label{ab8}
\|z^k-z^*\|_{G(\sigma)}^2 \geq \|z^{k+1}-z^*\|_{G(\sigma)}^2+\|z^{k+1}-z^{k+1/2}\|^2_{N(\sigma)}+\|z^{k+1/2}-z^k\|_{H(\sigma)}^2.
\end{equation}

Summing the inequality (\ref{ab8}) over $k$ from $0$ to $N$, we obtain
\begin{equation}\label{ab9}
\|z^0-z^*\|^2_{G(\sigma)}\geq \|z^{N+1}-z^*\|^2_{G(\sigma)}+\displaystyle
\sum_{k=1}^N\|z^{k+1}-z^{k+1/2}\|^2_{N(\sigma)}+\displaystyle
\sum_{k=1}^N\|z^{k+1/2}-z^{k}\|^2_{H(\sigma)}.
\end{equation}
Thus we obtain from (\ref{ab9}) that
$$
\begin{array}{l}
\|z^0-z^*\|^2_{G(\sigma)}\geq \|z^{N+1}-z^*\|^2_{G(\sigma)},\\[6pt]
\displaystyle
\sum_{k=1}^{\infty}\|z^{k+1}-z^{k+1/2}\|^2_{N(\sigma)}<\infty,\\[6pt]
\displaystyle
\sum_{k=1}^{\infty}\|z^{k+1/2}-z^{k}\|^2_{H(\sigma)}\leq \infty,
\end{array}
$$
 implying that
\begin{equation}\label{a10}
\|z^{k+1}-z^{k}\|\rightarrow 0.
\end{equation}
Since the sequence $\{z^{k}\}$ is bounded, there exist an element $\bar z$ and $\{k_i \}\subset \textbf{N}$ such that
$z^{k_i} \rightarrow \bar z$. It follows from (\ref{a2}) that, for $z \in {\rm dom} f \times {\rm dom} g$,
$$
\displaystyle \frac{1}{2}\|z- z^k\|^2_{\Theta}+\sigma \widehat \Phi_L(z,z^k; z^k)-\displaystyle \frac{1}{2}\|z-z^{k+1/2}\|^2_{\widehat \Sigma+\Theta}\geq \sigma \widehat \Phi_L(z^{k+1/2}, z^k;z^k)+\displaystyle \frac{1}{2}
\|z^{k+1/2}-z^k\|^2_{\Theta},
$$
which is equivalent to
$$
\begin{array}{l}
\displaystyle \frac{1}{2}\|z- z^k\|^2_{\widehat \Sigma+\Theta}-\displaystyle \frac{1}{2}\|z-z^{k+1/2}\|^2_{\widehat \Sigma+\Theta}+\sigma [f(x)+K(z^k)+\langle {\rm D}_xK(z^k),x-x^k\rangle+\displaystyle \frac{1}{2}\|x- x^k\|^2_{\widehat \Sigma_f}-g(y^k)]
\\[6pt]
-\sigma[f(x^k)+K(z^k)+\langle {\rm D}_yK(z^k),y-y^k\rangle-\displaystyle \frac{1}{2}\|y- y^k\|^2_{\widehat \Sigma_g}-g(y)]
\geq \sigma \widehat \Phi_L(z^{k+1/2}, z^k;z^k)+\displaystyle \frac{1}{2} \|z^{k+1/2}-z^k\|^2_{\Theta}.
\end{array}
$$
The above relation indicates that
$$
\begin{array}{l}
\displaystyle \frac{1}{2}\|z- z^k\|^2_{\sigma\widehat \Sigma+\Theta}-\displaystyle \frac{1}{2}\|z-z^{k+1/2}\|^2_{\widehat \Sigma+\Theta}+\sigma [f(x)+\langle {\rm D}_xK(z^k),x-x^k\rangle-g(y^k)]\\[6pt]
\quad -\sigma [f(x^k)+\langle {\rm D}_yK(z^k),y-y^k\rangle-g(y)]\\[6pt]
\geq \displaystyle \frac{1}{2}
\|z^{k+1/2}-z^k\|^2_{\Theta}+\sigma [f(x^{k+1/2})+\langle {\rm D}_xK(z^k),x^{k+1/2}-x^k\rangle+\displaystyle \frac{1}{2}\|x^{k+1/2}- x^k\|^2_{\widehat \Sigma_f}-g(y^k)]\\[6pt]
\quad -\sigma[f(x^k)+\langle {\rm D}_yK(z^k),y^{k+1/2}-y^k\rangle-\displaystyle \frac{1}{2}\|y^{k+1/2}- y^k\|^2_{\widehat \Sigma_g}-g(y^{k+1/2})].
\end{array}
$$
Thus, as $\sigma \in (0,1/2)$, we can get that
\begin{equation}\label{eq:b5}
\begin{array}{l}
\displaystyle \frac{1}{2}\|z- z^k\|^2_{\widehat \Sigma+\Theta}-\displaystyle \frac{1}{2}\|z-z^{k+1/2}\|^2_{\widehat \Sigma+\Theta}+\sigma [f(x)+\langle {\rm D}_xK(z^k),x-x^k\rangle]\\[12pt]
\quad  -\sigma [\langle {\rm D}_yK(z^k),y-y^k\rangle-g(y)]\\[6pt]
\geq \displaystyle \frac{1}{2}
\|z^{k+1/2}-z^k\|^2_{\Theta}+\sigma [\langle {\rm D}_xK(z^k),x^{k+1/2}-x^k\rangle+\displaystyle \frac{1}{2}\|x^{k+1/2}- x^k\|^2_{\widehat \Sigma_f}]\\[6pt]
\quad -\sigma[\langle {\rm D}_yK(z^k),y^{k+1/2}-y^k\rangle-\displaystyle \frac{1}{2}\|y^{k+1/2}- y^k\|^2_{\widehat \Sigma_g}]+\sigma f(x^{k+1/2})+\sigma g(y^{k+1/2}).\\[6pt]
\end{array}
\end{equation}
Since $f$ and $g$ are lower semi-continuous and $\|z^{k_i+1/2}-z^{k_i}\|\rightarrow 0$ as $i \rightarrow \infty$, taking the lower limit along $\{k_i\}$ on both sides of (\ref{eq:b5}),  we obtain
\begin{equation}\label{eq:b6}
\begin{array}{l}
 [f(x)+\langle {\rm D}_xK(\bar z),x-\bar x\rangle-g(\bar y)]
-[f(\bar x)+\langle {\rm D}_yK(\bar z),y-\bar y\rangle-g(y)]
\geq 0.
\end{array}
\end{equation}
In view of the convexity-concavity of $K$, we have that
$$
\langle {\rm D}_xK(\bar z),x-\bar x\rangle \leq K(x, \bar y)-K(\bar x,\bar y)
$$
and
$$
\langle {\rm D}_yK(\bar z),y-\bar y\rangle \geq K(\bar x,y)-K(\bar x,\bar y).
$$
Combining these with (\ref{eq:b6}), we obtain
$$
[f(x)+ K(x, \bar y)-g(\bar y)]
-[f(\bar x)+K(\bar x,y)-g(y)]
\geq 0.
$$
This implies that $\bar z \in \overline \Omega$.

Thus, any limit
point of the sequence $\{z^k\}$
is a solution of the problem. For any limit point $z^*$ of $\{z^k\}$, the relation (\ref{ab8}) indicates that the quantity $\|z^k-z^*\|^2_{G(\sigma)}$ decreases
monotonically. Combing these two facts, we know that the sequence $\{z^k\}$
can have only one limit
point, that is, $z^k$
converges monotonically with respect to $G(\sigma)$-norm to one of the solutions of the
problem; {\it i.e.}, $z^k \rightarrow z^*$. This proves the statement.\hfill $\Box$
\begin{remark}\label{rem:3.1} Although some assumptions on $\sigma$, $\Theta$ and $\widehat \Sigma$ are required in Theorem \ref{gconv}, they are reasonable.
\begin{itemize}
\item[a)]
In (\ref{eq:par}), we require an condition about the range of $\sigma$:
$$0< \sigma <\displaystyle \min\left\{ \displaystyle \frac{\lambda_{\min}(\widehat \Sigma+\Theta)}{\sqrt{2}\widehat\eta_0},\displaystyle \frac{1}{2}\right\},$$
this is not a strict condition, because $\sigma^{-1}$ plays an role as a penalty and usually $\sigma$ is a small positive real number.
\item [b)]In (\ref{eq:par}), we require an condition about $\Theta$:
$$
 \Theta \succ (\widehat \eta_0+2)\sigma I.
 $$
 This is also a reasonable condition, because $\Theta$ is an artificial matrix, which is provided by the  user. Certainly, the user may choose a positively definite operator and this condition can easily be satisfied if $\sigma>0$ small enough.
 \item [c)]
In (\ref{eq:hatsigma}), $\widehat \Sigma$ is required to satisfy
$$
\widehat \Sigma \succ \displaystyle \frac{6\eta_0\sigma}{1-2\sigma} I,
$$
this is not a restricted assumption because we usually choose a large positively definite operator $\widehat \Sigma$ so that
it satisfies (\ref{def-LB-f}) and (\ref{def-LB-g}).
\end{itemize}
\end{remark}
\section{The rate of convergence}
\setcounter{equation}{0}
\quad \, Under  Assumption \ref{K-cc}, $L$ is a convex-concave function and $(\bar x,\bar y) \in \overline \Omega$ if and only if $(\bar x,\bar y)$ satisfies
\begin{equation}\label{eq:KKT}
\left
\{
\begin{array}{l}
0 \in \partial f(x)+{\rm D}_x K(x,y),\\[6pt]
0\in \partial g(y)-{\rm D}_y  K(x,y).
\end{array}
\right.
\end{equation}
This is a generalized equation version for the optimality. We can also express the optimality as an equation
$$
R(z)=R(x,y)=0,
$$
where
\begin{equation}\label{eq:ST}
R(z)=\left
[
\begin{array}{l}
x-{\bf P}_f (x-{\rm D}_x K(x,y))\\[4pt]
y-{\bf P}_g (y+{\rm D}_y K(x,y))
\end{array}
\right
]
\end{equation}
and ${\bf P}_{\psi}$ is the proximal mapping. Here, for a convex function $\psi$, ${\bf P}_{\psi}$ is defined by
$$
{\bf P}_{\psi}(w)=\arg\min_{w'}\left\{\psi(w')-\displaystyle \frac{1}{2}\|w'-w\|^2\right\}.
$$
Then we can express the set of saddle points of $L$ as
$$
\overline \Omega = R^{-1}(0).
$$
To develop the rate of convergence of Algorithm \ref{alg-1}, we need the following metric subregularity of $R$ at $(z^*,0)$.
\begin{assumption}\label{ass-MSR}
Suppose that $R$  is locally metrically subregular at $(z^*,0)$; {\it i.e.}, there exist  $\varepsilon_0>0$ and $\kappa_0>0$ such that
\begin{equation}\label{eq:MSR}
{\rm dist}\,(z, R^{-1}(0)=\overline \Omega) \leq \kappa_0 \|R(z)\|, \ \forall z \in \textbf{B}(z^*,\varepsilon_0).
\end{equation}
\end{assumption}
\begin{proposition}\label{prop:RE}
Let Assumption \ref{assu-K} be satisfied.  Let $\{z^k=(x^k,y^k)\}$ and $\{z^{k+1/2}=(x^{k+1/2},y^{k+1/2})\}$  be generated by Algorithm
\ref{alg-1}. Then for $k \geq 0$,
\begin{equation}\label{eq:REs}
\|R(z^{k+1})\|^2\leq \|z^{k+1}-z^{k+1/2}\|^2_{6\eta_0^2 I+3(\widehat \Sigma+\Theta)^*(\widehat \Sigma+\Theta)}+\|z^{k+1/2}-z^{k}\|^2_{3\Theta^*\Theta}.
\end{equation}
\end{proposition}
{\bf Proof}. From the definition of $z^{k+1}$ in (\ref{key-o}) and the optimality condition, we obtain
$$
\begin{array}{l}
0 \in {\cal S}(x^{k+1}-x^k)+{\rm D}_x \widehat K(x^{k+1},y^{k+1/2};z^{k+1/2})+\partial f(x^{k+1}),\\[4pt]
0\in {\cal T}(y^{k+1}-y^k)-{\rm D}_y \widehat K(x^{k+1/2},y^{k+1};z^{k+1/2})+\partial g(y^{k+1}).
\end{array}
$$
Then by the definition of $\widehat K (z;z^{k+1/2})$, we know that
\begin{equation}\label{eq:c1}
\begin{array}{l}
0 \in {\cal S}(x^{k+1}-x^k)+{\rm D}_x K(z^{k+1/2})+
\widehat \Sigma_f(x^{k+1}-x^{k+1/2})+\partial f(x^{k+1}),\\[4pt]
0\in {\cal T}(y^{k+1}-y^k)-{\rm D}_y  K(z^{k+1/2})+\widehat \Sigma_g(y^{k+1}-y^{k+1/2})+\partial g(y^{k+1}),
\end{array}
\end{equation}

$$
\begin{array}{l}
x^{k+1}={\bf P}_f (x^{k+1}-{\rm D}_x K(z^{k+1/2})-\widehat \Sigma_f(x^{k+1}-x^{k+1/2})- {\cal S}(x^{k+1}-x^k)),\\[4pt]
y^{k+1}={\bf P}_g (y^{k+1}+{\rm D}_y K(z^{k+1/2})-\widehat \Sigma_g(y^{k+1}-y^{k+1/2})- {\cal T}(y^{k+1}-y^k)).
\end{array}
$$
Therefore we obtain
$$
\begin{array}{lll}
\|R(z^{k+1})\|^2 & \leq & \|{\rm D}_x K(x^{k+1},y^{k+1})-{\rm D}_x K(z^{k+1/2})-\widehat \Sigma_f(x^{k+1}-x^{k+1/2})- {\cal S}(x^{k+1}-x^k)\|^2\\[6pt]
&&+ \|-{\rm D}_y K(x^{k+1},y^{k+1})+{\rm D}_y K(z^{k+1/2})-\widehat \Sigma_g(y^{k+1}-y^{k+1/2})-{\cal T}(y^{k+1}-y^k)\|^2\\[6pt]
&=& \|[{\rm D}_x K(z^{k+1})-{\rm D}_x K(z^{k+1/2})]- ({\cal S}+\widehat \Sigma_f)(x^{k+1}-x^{k+1/2})-{\cal S}(x^{k+1/2}-x^k)\|^2\\[6pt]
&&+ \|[{\rm D}_y K(z^{k+1/2})- {\rm D}_y K(z^{k+1}]-({\cal T}+\widehat \Sigma_g)(y^{k+1}-y^{k+1/2})-{\cal T}(y^{k+1/2}-y^k)\|^2\\[6pt]
&\leq & 3 \eta_0^2\|z^{k+1}-z^{k+1/2}\|^2+3\|x^{k+1}-x^{k+1/2}\|^2_{({\cal S}+\widehat \Sigma_f)^*({\cal S}+\widehat \Sigma_f)}+3\|x^{k+1/2}-x^{k}\|^2_{{\cal S}^*{\cal S}}\\[6pt]
&& +3 \eta_0^2\|z^{k+1}-z^{k+1/2}\|^2+3\|y^{k+1}-y^{k+1/2}\|^2_{({\cal T}+\widehat \Sigma_g)^*({\cal T}+\widehat \Sigma_g)}+3\|y^{k+1/2}-y^{k}\|^2_{{\cal T}^*{\cal T}}\\[6pt]
&=& \|z^{k+1}-z^{k+1/2}\|^2_{6\eta_0^2 I+3(\widehat \Sigma+\Theta)^*(\widehat \Sigma+\Theta)}+\|z^{k+1/2}-z^{k}\|^2_{3\Theta^*\Theta},
\end{array}
$$
which implies the truth of the statement.  \hfill $\Box$\\

Now we define
$$
{\rm dist}_{\widehat \Sigma+\Theta}\, (z,\overline \Omega)=\inf_{z' \in \overline \Omega} \{\|z'-z\|_{\widehat \Sigma+\Theta}\}.$$
With the help of Proposition \ref{prop:RE}, we can prove the linear rate of convergence under the locally metrical subregularity of
 $R$  at $(z^*,0)$.
\begin{theorem}\label{gconR}
Let Assumptions \ref{K-cc}, \ref{assu-K}, \ref{ass-o} and  \ref{ass-MSR} be satisfied at every point $z^* \in \overline \Omega$.  Let $\{z^k=(x^k,y^k)\}$ and $\{z^{k+1/2}=(x^{k+1/2},y^{k+1/2})\}$  be generated by Algorithm
\ref{alg-1}.
Suppose that ${\cal S}$ and ${\cal T}$ satisfy (\ref{cond-2}) or equivalently $\Theta \succeq 0$ and $\Theta+\widehat \Sigma\succ 0$.
Suppose also that $\sigma$, ${\cal S}$ and ${\cal T}$ satisfy
\begin{equation}\label{eq:para}
0< \sigma <\displaystyle \min\left\{ \displaystyle \frac{\lambda_{\min}(\widehat \Sigma+\Theta)}{\sqrt{2}\widehat\eta_0},\displaystyle \frac{1}{2}\right\},\, \Theta \succ (\widehat \eta_0+2)\sigma I.
\end{equation}
If $\widehat \Sigma$ satisfies
\begin{equation}\label{eq:hatsigmaa}
\widehat \Sigma \succ \displaystyle \frac{6\eta_0\sigma}{1-2\sigma} I,
\end{equation}
then $\{z^k=(x^k,y^k)\}$ converges linearly with respect to some norm to an element of $\overline \Omega$.
\end{theorem}
{\bf Proof}. From Theorem \ref{gconv}, $\{z^k\}$ converges to an element of  $\overline \Omega$, say $\overline z$. This indicates that
$$
z^k \in  \textbf{B}(\overline z,\varepsilon_0)\subset \textbf{B}(\overline \Omega,\varepsilon_0)
$$
for any $k> N$, where $N$ is a large integer.

Recall from the proof of Theorem \ref{gconv} that $G(\sigma)$,$N(\sigma)$ and $H(\sigma)$ are positively definite, where
$$
\begin{array}{l}
G(\sigma)=\widehat \Sigma+\Theta-6\eta_0\sigma I,\\[4pt]
N(\sigma)=\widehat \Sigma+\Theta-(\widehat \eta_0+6\eta_0)\sigma I,\\[4pt]
H(\sigma)=\Theta-(\widehat \eta_0+2)\sigma I.
\end{array}
$$
Also recall from (\ref{ab8}) that we have the following relation
\begin{equation}\label{ab8a}
\|z^k-z^*\|_{G(\sigma)}^2 \geq \|z^{k+1}-z^*\|_{G(\sigma)}^2+\|z^{k+1}-z^{k+1/2}\|^2_{N(\sigma)}+\|z^{k+1/2}-z^k\|_{H(\sigma)}^2.
\end{equation}
 Since both $N(\sigma)$ and $H(\sigma)$ are positively definite,  there must exist a positive number $\mu>0$ such that
$$
6\eta_0^2 I+3(\widehat \Sigma+\Theta)^*(\widehat \Sigma+\Theta)\prec \mu  N(\sigma),\quad  3\Theta^*\Theta\prec \mu H(\sigma).
$$
Then for $k >N$, from Assumption \ref{ass-MSR}, we have from Proposition \ref{prop:RE} for any $z^* \in \overline \Omega$ that
\begin{equation}\label{ab1a}
\begin{array}{rcl}
{\rm dist}_{G(\sigma)}\, (z^{k+1},\overline \Omega)^2 &
\leq & \lambda_{\max}(G(\sigma)) {\rm dist}(z^{k+1},\overline \Omega)^2\\[10pt]
& \leq & \kappa_0^2 \lambda_{\max}(G(\sigma))\|R(z^{k+1}\|^2\\[10pt]
& \leq & \kappa_0^2\lambda_{\max}(G(\sigma))\left[\|z^{k+1}-z^{k+1/2}\|^2_{6\eta_0^2 I+3(\widehat \Sigma+\Theta)^*(\widehat \Sigma+\Theta)}+\|z^{k+1/2}-z^{k}\|^2_{3\Theta^*\Theta}\right]\\[12pt]
&\leq & \kappa_0^2\mu\lambda_{\max}(G(\sigma))\left[ \|z^{k+1}-z^{k+1/2}\|^2_{N(\sigma)}+\|z^{k+1/2}-z^{k}\|^2_{H(\sigma)}\right]\\[6pt]
& \leq & \kappa_0^2\mu\lambda_{\max}(G(\sigma))\left[\|z^k-z^*\|_{G(\sigma)}^2-\|z^{k+1}-z^*\|_{G(\sigma)}^2\right],
\end{array}
\end{equation}
where the last inequality comes from (\ref{ab8a}).
Taking $z^*$ in (\ref{ab1a}) as
$$
z^* =\arg\min\left\{\|z-z^k\|_{G(\sigma)}:z \in \overline \Omega\right\},
$$
and noting that
$$
{\rm dist}_{G(\sigma)}\, (z^{k+1},\overline \Omega)^2\leq \|z^{k+1}-z^*\|_{G(\sigma)}^2,
$$
we get from (\ref{ab1a}) that
$$
{\rm dist}_{G(\sigma)}\, (z^{k+1},\overline \Omega)^2+ \kappa_0^2\mu\lambda_{\max}(G(\sigma)) \|z^{k+1}-z^*\|_{G(\sigma)}^2\leq \kappa_0^2\mu \lambda_{\max}(G(\sigma)) {\rm dist}_{G(\sigma)}\, (z^k,\overline \Omega)^2.
$$
The above relation indicates that
$$
{\rm dist}_{G(\sigma)}\, (z^{k+1},\overline \Omega)\leq \displaystyle \frac{1}{\sqrt{1+[\kappa_0^2\mu\lambda_{\max}(G(\sigma))]^{-1} }}
{\rm dist}_{G(\sigma)}\, (z^k,\overline \Omega).
$$
This means that $z^k$ converges to an element of $\overline \Omega$ with linear rate of convergence with respect to $G(\sigma) $-norm. The proof is completed. \hfill $\Box$\\

\section{Numerical Experiments}\label{Sec6}
%
In this section, we present some preliminary numerical experiments to illustrate the performance of Algorithm \ref{alg-1}. All numerical experiments are implemented by MATLAB R2019a on a laptop with Intel(R) Core(TM) i5-6200U 2.30GHz and 8GB memory. We discuss the application of the majorized semi-proximal alternating coordinate method in four different forms of minimax optimization problems. The advantages of mspACM are shown  for two cases when $f$ and $g$ are smooth or nonsmooth.

\subsection{Smooth Saddle Point Problems}
When $f$ and $g$ are smooth, the minimax optimization problem can be abbreviated as the following form
\begin{equation}\label{sspp}
\min_{x \in \Re^n}\max_{y \in \Re^m}\,K(x,y).
\end{equation}
Many methods can solve saddle point problems, such as the proximal point (PP) method (see \cite{Roc1976b}), the optimistic gradient descent ascent (OGDA) method and the extra-gradient (EG) method (see \cite{Ozd2019a}). We compare the proposed method, mspACM, with these methods for  the linear regression problem. For solving Problem (\ref{sspp}), the updating formula of mspACM can be simplified as
\begin{equation}\label{key-o1}
\begin{array}{ccl}
z^{k+1/2}&=&-\sigma(\sigma\widehat\Sigma+\Theta)^{-1}\widetilde{{\rm D}} K(z^k)+z^k,\\[4pt]
z^{k+1}&=&(\sigma\widehat\Sigma+\Theta)^{-1}[-\widetilde{{\rm D}} K(z^{k+1/2})+\sigma\widehat\Sigma (z^{k+1/2})+\Theta (z^k)],
\end{array}
\end{equation}
where $z^k=(x^{k},y^{k})$ and $\widehat\Sigma$ and $\Theta$ are two linear operators defined in \eqref{hatSigma} and \eqref{Theta}, respectively. The operator $\widetilde{\rm D} K(\cdot)$ is expressed as
$$\widetilde{\rm D} K(z)=\left(
             \begin{array}{c}
               {\rm D}_x K(x,y) \\
               - {\rm D}_y K(x,y)\\
             \end{array}
           \right)
$$
for any $z=(x,y)\in\Re^n\times\Re^m $.

The saddle point reformulation of the linear regression is of the form
\begin{equation}\label{lr}
\min_{x \in \Re^n}\max_{y \in \Re^m}\,\frac{1}{m}\left[-\frac{1}{2}\|y\|^2-b^Ty+y^TAx\right]+\frac{\lambda}{2}\|x\|^2.
\end{equation}
We set $n=m$ and the rows of the matrix $A$  are generated by a Gaussian distribution $\mathcal {N}(0,I_n)$. Let $b=0$ and $\lambda=1/m$. In PP, EG and OGDA, parameters and step sizes are selected for best performance. In mspACM, we set $\sigma=1$. The linear operators are chosen as ${\cal S}={\cal T}=\|A\|_2I_n$, $\widehat \Sigma_f =0.1A^TA$ and $\widehat \Sigma_g =0.1AA^T$.

\begin{figure}[htbp]
 \centering

 \subfigure[n=10]{
 \begin{minipage}[t]{0.5\linewidth}
 \centering
 \includegraphics[width=1\textwidth]{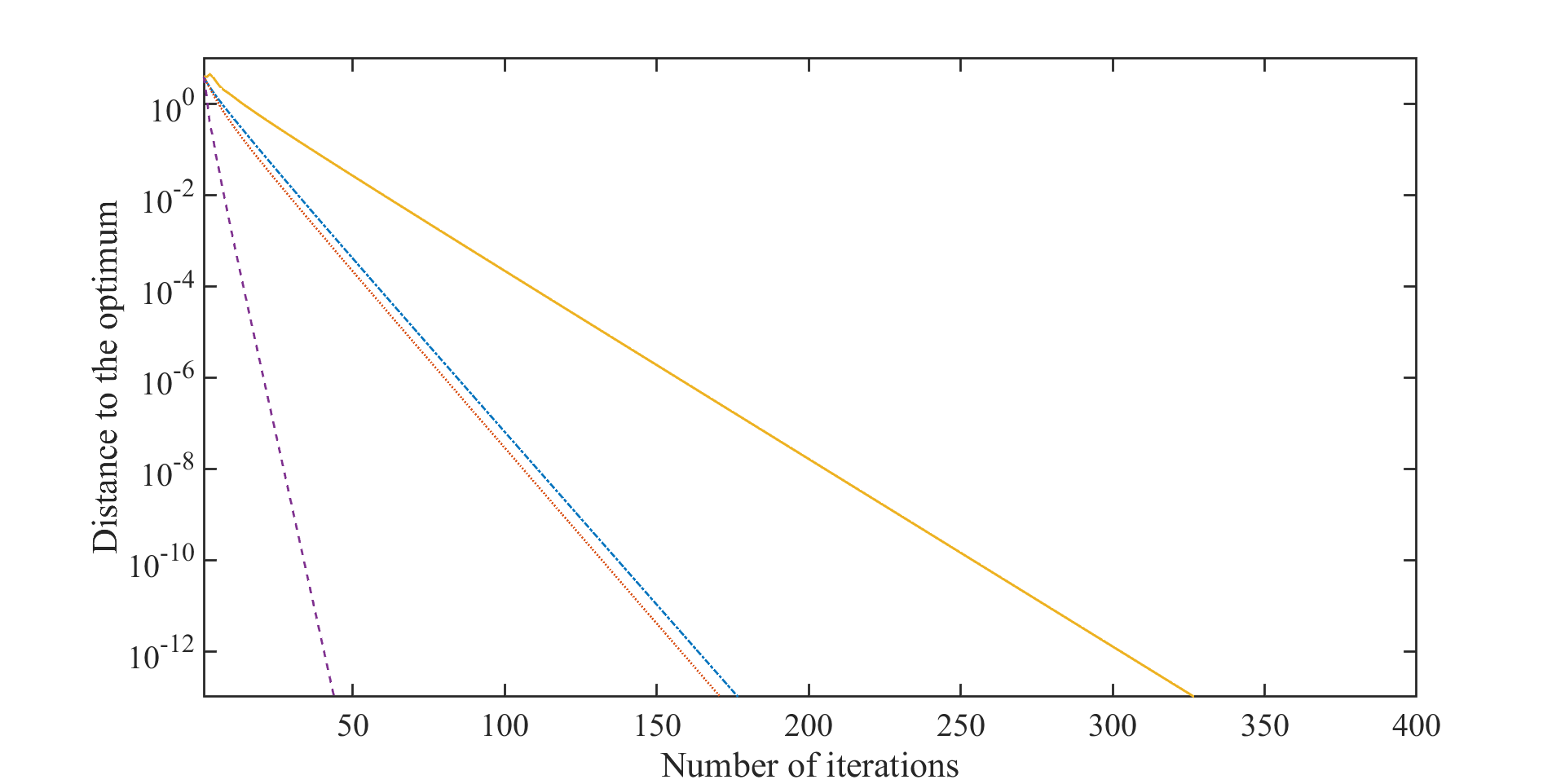}
 \end{minipage}%
 }%
\subfigure[n=100]{
 \begin{minipage}[t]{0.5\linewidth}
 \centering
 \includegraphics[width=1\textwidth]{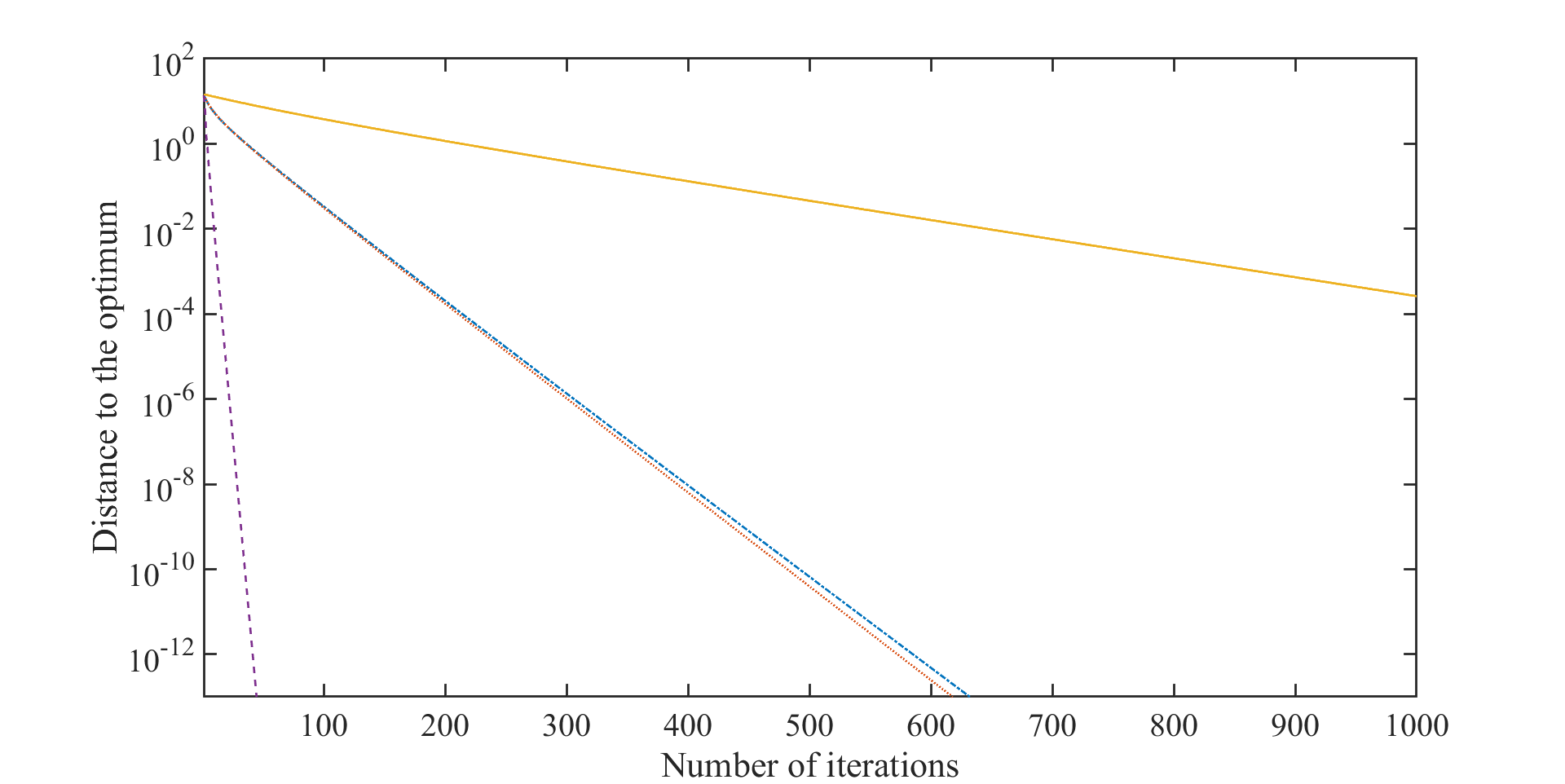}
 \end{minipage}%
 }%

 \subfigure[n=1000]{
 \begin{minipage}[t]{0.5\linewidth}
 \centering
 \includegraphics[width=1\textwidth]{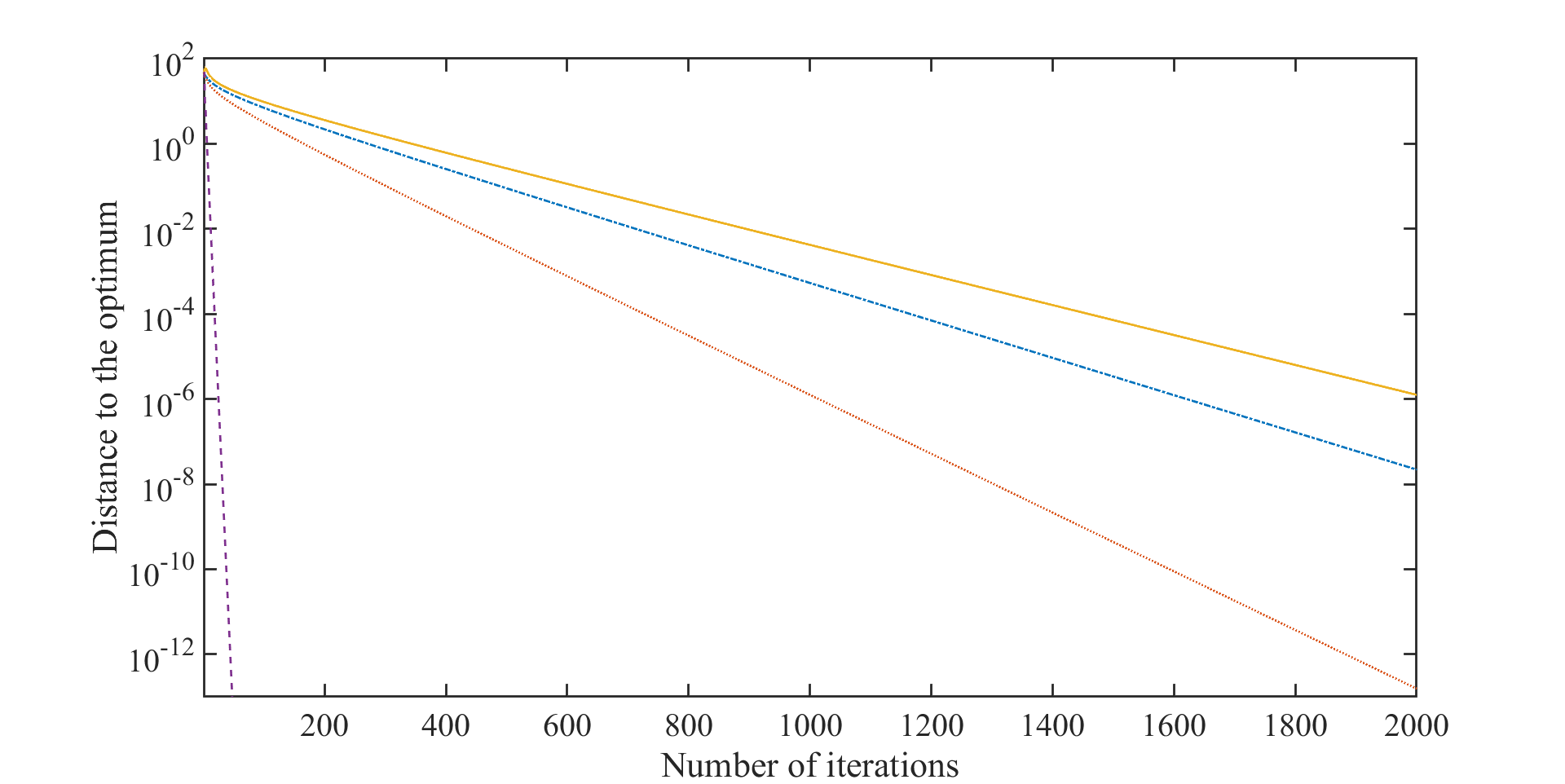}
 \end{minipage}%
 }%

 \centering
 \caption{ Compare mspACM, PP, EG and OGDA in terms of number of iterations under different dimensions for the linear regression. Step sizes of EG and OGDA are tuned for best performance.}\label{fig1}
 \end{figure}

In Figure \ref{fig1}, we compare the performances of the four methods with respect to the number of iterations when the dimension varies from $n=10,\ 100$ and $1000$. The same initial point $x^0$ is chosen. Generally speaking, all the four methods converge linearly to the optimal solution, and the proximal point (PP) method has the best performance. Our method, mspACM, is the second best, which converges faster than EG and OGDA. It can be observed that for the low-dimensional strongly convex-strongly concave saddle point problem \eqref{lr}, the convergence rate of EG is very close to mspACM. When the dimension $n$ increases, the performance  of mspACM becomes much better then  EG and OGDA.

It  is reasonable to explain the best performance of PP because the proximal point method is asymptotically superlinear and it has an explicit solution for every subproblem when solving this simple problem.

\subsection{Nonsmooth Convex Minimax Optimization Problems}
In this part, we focus on the convex minimax optimization problem of the form
\begin{equation}\label{cmm}
\min_{x \in X}\max_{y \in Y}\,\mu_x\|x\|_{\infty}+\frac{\lambda}{2}\|x\|^2+\frac{1}{m}\left[-\frac{1}{2}\|y\|^2-b^Ty+y^TAx\right]-\mu_y\|y\|_{\infty},
\end{equation}
where $X\in \Re^n$ and $Y\in \Re^m$ are two convex sets and $\|\cdot\|_{\infty}$ represents the infinite norm of a finite-dimensional vector space.

We can not use EG and OGDA to solve problem (\ref{cmm}) as this is a nonsmooth minimax optimization problem.  It is quite difficult to use the proximal method to solve problem (\ref{cmm}) because its  subproblems can not be solved directly. However, although mspACM can not provide the explicit solutions of the subproblems, they are simple nonsmooth unconstrained optimization problems and can be solved by Matlab code ``fminunc" efficiently.

\emph{Unconstrained convex minimax optimization problem}. We set $n=m$, $b=0$, $\lambda=1/m$ and $\mu_x=\mu_y=1$. In this case, the linear operator ${\cal S}={\cal T}=\|A\|_2I_n$ and $\widehat \Sigma_f =0.1A^TA,\ \widehat \Sigma_g =0.1AA^T$. The convergence of mspACM is shown in Table \ref{tab1}, with respect to the number of iterations and CPU time under different dimensions, different condition numbers of matrix $A$, and different values of the parameter $\sigma$. It is easy to see that the convergence rate becomes  significantly slower as the dimension $n$ increases. However, the condition number of matrix $A$ has very little effect on the convergence of mspACM. On the other hand, within the same CPU time, the convergence of mspACM becomes faster as $\sigma$ increases.

\begin{sidewaystable}[h]
\centering
\begin{tabular}{ccc|c|c|c|c|c|c|c|c|c|c|}
\cline{4-13}
                                                &                                                      &              & \multicolumn{2}{c|}{$\epsilon=10^{-1}$} & \multicolumn{2}{c|}{$\epsilon=10^{-3}$} & \multicolumn{2}{c|}{$\epsilon=10^{-5}$} & \multicolumn{2}{c|}{$\epsilon=10^{-7}$} & \multicolumn{2}{c|}{$\epsilon=10^{-9}$} \\ \cline{4-13}
                                                &                                                      &              & $T$                 & $t$(sec)                & $T$                & $t$(sec)                 & $T$                & $t$(sec)                 & $T$                & $t$(sec)                & $T$     & $t$(sec)                            \\ \hline
\multicolumn{1}{|c|}{}                          & \multicolumn{1}{c|}{}                                & $\sigma=1$   & 9                 & 6.80                & 11               & 7.56                 & 12               & 7.95                 & 13               & 8.19                 & 14    & 8.37                            \\ \cline{3-13}
\multicolumn{1}{|c|}{}                          & \multicolumn{1}{c|}{\multirow{-2}{*}{$\kappa=10$}}   & $\sigma=0.1$ & 77                & 15.27               & 100              & 23.73                & 104              & 24.91                & 105              & 25.17                & 106   & 25.30                           \\ \cline{2-13}
\multicolumn{1}{|c|}{}                          & \multicolumn{1}{c|}{}                                & $\sigma=1$   & 5                 & 2.93                & 7                & 3.82                 & 8                & 4.60                 & 8                & 4.60                 & 9     & 5.22                            \\ \cline{3-13}
\multicolumn{1}{|c|}{}                          & \multicolumn{1}{c|}{\multirow{-2}{*}{$\kappa=50$}}   & $\sigma=0.1$ & 46                & 10.34               & 61               & 14.94                & 63               & 15.47                & 64               & 15.74                & 65    & 16.01                           \\ \cline{2-13}
\multicolumn{1}{|c|}{}                          & \multicolumn{1}{c|}{}                                & $\sigma=1$   & 5                 & 3.27                & 6                & 3.81                 & 7                & 4.41                 & 8                & 5.00                 & 9     & 5.46                            \\ \cline{3-13}
\multicolumn{1}{|c|}{\multirow{-6}{*}{$n=10$}}  & \multicolumn{1}{c|}{\multirow{-2}{*}{$\kappa=200$}}  & $\sigma=0.1$ & 40                & 9.69                & 51               & 12.43                & 53               & 12.71                & 54               & 12.92                & 55    & 13.08                           \\ \hline
\multicolumn{1}{|c|}{}                          & \multicolumn{1}{c|}{}                                & $\sigma=1$   & 13                & 19.95               & 26               & 40.16                & 32               & 48.56                & 38               & 53.81                & 44    & 58.86                           \\ \cline{3-13}
\multicolumn{1}{|c|}{}                          & \multicolumn{1}{c|}{\multirow{-2}{*}{$\kappa=100$}}  & $\sigma=0.1$ & 19                & 20.42               & 45               & 75.81                & 54               & 86.48                & 60               & 92.86                & 65    & 97.01                           \\ \cline{2-13}
\multicolumn{1}{|c|}{}                          & \multicolumn{1}{c|}{}                                & $\sigma=1$   & 14                & 21.72               & 27               & 42.52                & 33               & 51.84                & 39               & 58.14                & 43    & 61.68                           \\ \cline{3-13}
\multicolumn{1}{|c|}{}                          & \multicolumn{1}{c|}{\multirow{-2}{*}{$\kappa=1000$}} & $\sigma=0.1$ & 19                & 16.52               & 43               & 57.88                & 49               & 67.90                & 51               & 70.26                & 55    & 73.53                           \\ \cline{2-13}
\multicolumn{1}{|c|}{}                          & \multicolumn{1}{c|}{}                                & $\sigma=1$   & 16                & 24.15               & 31               & 48.29                & 35               & 53.20                & 41               & 59.80                & 43    & 61.26                           \\ \cline{3-13}
\multicolumn{1}{|c|}{\multirow{-6}{*}{$n=50$}}  & \multicolumn{1}{c|}{\multirow{-2}{*}{$\kappa=5000$}} & $\sigma=0.1$ & 19                & 19.71               & 43               & 58.79                & 52               & 73.34                & 59               & 81.20                & 63    & 83.97                           \\ \hline
\multicolumn{1}{|c|}{}                          & \multicolumn{1}{c|}{}                                & $\sigma=1$   & 29                & 33.11               & 75               & 157.41               & 100              & 235.41               & 112              & 267.93               & 116   & 273.04                          \\ \cline{3-13}
\multicolumn{1}{|c|}{}                          & \multicolumn{1}{c|}{\multirow{-2}{*}{$\kappa=100$}}  & $\sigma=0.1$ & 335               & 168.64              & 750              & 1057.81              & 882              & 1509.52              & 890              & 1524.81              & 894   & 1530.36                         \\ \cline{2-13}
\multicolumn{1}{|c|}{}                          & \multicolumn{1}{c|}{}                                & $\sigma=1$   & 29                & 38.48               & 79               & 181.79               & 108              & 278.94               & 115              & 294.65               & 119   & 300.58                          \\ \cline{3-13}
\multicolumn{1}{|c|}{}                          & \multicolumn{1}{c|}{\multirow{-2}{*}{$\kappa=1000$}} & $\sigma=0.1$ & 351               & 158.66              & 735              & 1043.61              & 867              & 1478.71              & 876              & 1496.09              & 880   & 1502.40                         \\ \cline{2-13}
\multicolumn{1}{|c|}{}                          & \multicolumn{1}{c|}{}                                & $\sigma=1$   & 28                & 39.89               & 78               & 187.14               & 111              & 292.06               & 120              & 311.01               & 124   & 317.16                          \\ \cline{3-13}
\multicolumn{1}{|c|}{\multirow{-6}{*}{$n=100$}} & \multicolumn{1}{c|}{\multirow{-2}{*}{$\kappa=1000$}} & $\sigma=0.1$ & 360               & 184.59              & 752              & 1127.88              & 896              & 1582.17              & 908              & 1605.02              & 911   & 1609.378                        \\ \hline
\multicolumn{1}{|c|}{}                          & \multicolumn{1}{c|}{}                                & $\sigma=1$   & 40                & 436.56              & 68               & 742.15               & 82               & 887.87               & 92               & 963.47               & 95    & 976.87                          \\ \cline{3-13}
\multicolumn{1}{|c|}{}                          & \multicolumn{1}{c|}{\multirow{-2}{*}{$\kappa=100$}}  & $\sigma=0.1$ & 21                & 132.58              & 59               & 552.76               & 101              & 1029.34              & 114              & 1138.94              & 118   & 1154.44                         \\ \cline{2-13}
\multicolumn{1}{|c|}{}                          & \multicolumn{1}{c|}{}                                & $\sigma=1$   & 47                & 457.32              & 82               & 822.88               & 94               & 929.33               & 107              & 1009.06              & 111   & {\color[HTML]{000000} 1026.05}  \\ \cline{3-13}
\multicolumn{1}{|c|}{}                          & \multicolumn{1}{c|}{\multirow{-2}{*}{$\kappa=1000$}} & $\sigma=0.1$ & 22                & 152.49              & 60               & 566.99               & 103              & 1041.50              & 116              & 1146.41              & 121   & 1165.42                         \\ \cline{2-13}
\multicolumn{1}{|c|}{}                          & \multicolumn{1}{c|}{}                                & $\sigma=1$   & 44                & 415.08              & 65               & 634.97               & 80               & 791.97               & 91               & 871.00               & 95    & 888.81                          \\ \cline{3-13}
\multicolumn{1}{|c|}{\multirow{-6}{*}{$n=200$}} & \multicolumn{1}{c|}{\multirow{-2}{*}{$\kappa=10^5$}} & $\sigma=0.1$ & 21                & 149.81              & 59               & 560.90               & 106              & 1075.58              & 120              & 1173.37              & 125   & 1192.25                         \\ \hline
\end{tabular}
\caption{ Numerical results of mspACM for the minimax  problem \eqref{cmm}. $n=m={\rm dim}\,x$; $\kappa$-condition number of $A$; $\sigma$-the parameter of mspACM; $\epsilon$-the relative error (i.e., $\|(x^k;y^k)\|_2/\|(x^0;y^0)\|_2$); $T$- number of iterations; $t$-CPU time.}\label{tab1}
\end{sidewaystable}

\emph{Convex minimax optimization problem with convex constraints}. We tested three cases as follows.
\begin{description}
  \item[Case 1.] Linear equality constraints:
  $$
X=\{x\in \Re^n: B_xx=b_{ex}\},\quad Y=\{y\in {\Re^m}: B_yy=b_{ey}\}.
$$
  \item[Case 2.] Linear inequality constraints:
  $$
X=\{x\in \Re^n: M_xx\leq b_{ix}\},\quad Y=\{y\in {\Re^m}: M_yy\leq b_{iy}\}.
$$
  \item[Case 3.] Quadratic constraints:
  $$
X=\{x\in \Re^n: x^TQ_xx+c_x^Tx+b_{qx}\leq0\},\quad Y=\{y\in {\Re^m}: y^TQ_yy+c_y^Ty+b_{qy}\leq0\},
$$
where $Q_x$ and $Q_y$ are positive semi-definite matrices .
\end{description}
We rewrite Problem \eqref{cmm}  as
\begin{equation}\label{example-alt}
\min_{x \in {\Re^n}}\max_{y \in { \Re^m}}\,\delta_X(x)+\mu_x\|x\|_{\infty}+\frac{\lambda}{2}\|x\|^2+\frac{1}{m}\left[-\frac{1}{2}\|y\|^2-b^Ty+y^TAx\right]-\mu_y\|y\|_{\infty}+\delta_Y(y).
\end{equation}
For simplicity, let $b=0$, $\lambda=1/m$ and $\mu_x=\mu_y=1$. We consider the situation when $n\neq m$ in Problem \eqref{example-alt}. The linear operators are selected as ${\cal S}=\|A\|_2I_n,\ {\cal T}=\|A\|_2I_m$, $\widehat \Sigma_f =0.1A^TA$ and $\widehat \Sigma_g =0.1AA^T$ and the parameter $\sigma=1$.

The performance of mspACM under three different constraints is shown in Table \ref{tab2}. We can see from Table \ref{tab2} that, under each of the three different constraints, mspACM converges to the optimal solution of the Problem \eqref{example-alt} within a few number of iterations. Furthermore, it can be observed that as the dimensions $n$ and $m$ increase, the CPU time for implementing mspACM  increases  rapidly.

\begin{table}[h]
\centering
\begin{tabular}{cc|c|c|c|c|c|c|c|c|}
\cline{3-10}
                                                                                                 &                          & \multicolumn{2}{c|}{$\epsilon=10^{-1}$}    & \multicolumn{2}{c|}{$\epsilon=10^{-3}$}     & \multicolumn{2}{c|}{$\epsilon=10^{-5}$}     & \multicolumn{2}{c|}{$\epsilon=10^{-7}$}     \\ \cline{3-10}
                                                                                                 &                          & $T$                & $t$(sec)              & $T$                 & $t$(sec)              & $T$                 & $t$(sec)              & $T$                 & $t$(sec)              \\ \hline
\multicolumn{1}{|c|}{\multirow{3}{*}{\begin{tabular}[c]{@{}c@{}}$n=10$\\ $m=10$\end{tabular}}}   & Case 1. & 6 & 3.52 & 8  & 4.53 & 8  & 4.53 & 9  & 5.10                      \\ \cline{2-10}
\multicolumn{1}{|c|}{}                                                                           & Case 2. & 6 & 4.14 & 12 & 6.64 & 14 & 7.32 & 15 & 7.61 \\
 \cline{2-10}
\multicolumn{1}{|c|}{}                                                                           & Case 3. & 6 & 2.78 & 12 & 5.46 & 14 & 6.23 & 16 & 6.59 \\
 \hline
\multicolumn{1}{|c|}{\multirow{3}{*}{\begin{tabular}[c]{@{}c@{}}$n=20$\\ $m=50$\end{tabular}}}   & Case 1.                  & 4                  & 6.70                  & 5                   & 7.75                  & 6                   & 8.66                  & 7                   & 9.39                  \\ \cline{2-10}
\multicolumn{1}{|c|}{}                                                                           & Case 2.                  & 33                 & 35.98                 & 61                  & 63.88                 & 63                  & 65.93                 & 69                  & 69.65                 \\ \cline{2-10}
\multicolumn{1}{|c|}{}                                                                           & Case 3.                  & 45                 & 47.63                 & 80                  & 84.35                 & 86                  & 90.51                 & 92                  & 95.19                 \\ \hline
\multicolumn{1}{|c|}{\multirow{3}{*}{\begin{tabular}[c]{@{}c@{}}$n=50$\\ $m=20$\end{tabular}}}   & Case 1.                  & 5                  & 6.53                  & 6                   & 7.31                  & 7                   & 8.76                  & 8                   & 10.01                 \\ \cline{2-10}
\multicolumn{1}{|c|}{}                                                                           & Case 2.                  & 44                 & 45.11                 & 81                  & 81.90                 & 86                  & 86.78                 & 90                  & 90.57                 \\ \cline{2-10}
\multicolumn{1}{|c|}{}                                                                           & Case 3.                  & 19                 & 29.10                 & 26                  & 35.70                 & 29                  & 39.04                 & 34                  & 42.88                 \\ \hline
\multicolumn{1}{|c|}{\multirow{3}{*}{\begin{tabular}[c]{@{}c@{}}$n=100$\\ $m=100$\end{tabular}}} & Case 1.                  & 17                 & 51.72                 & 39                  & 121.54                & 46                  & 144.36                & 57                  & 165.73                \\ \cline{2-10}
\multicolumn{1}{|c|}{}                                                                           & Case 2.                  & 21                 & 65.15                 & 44                  & 140.91                & 52                  & 162.85                & 61                  & 179.35                \\ \cline{2-10}
\multicolumn{1}{|c|}{}                                                                           & Case 3.                  & 20                 & 62.80                 & 38                  & 123.08                & 45                  & 146.17                & 57                  & 173.13                \\ \hline
\multicolumn{1}{|c|}{\multirow{3}{*}{\begin{tabular}[c]{@{}c@{}}$n=100$\\ $m=300$\end{tabular}}} & Case 1.                  & 13                 & 213.40                & 48                  & 853.78                & 56                  & 996.23                & 66                  & 1089.59               \\ \cline{2-10}
\multicolumn{1}{|c|}{}                                                                           & Case 2.                  & 14                 & 236.87                & 30                  & 534.90                & 47                  & 836.47                & 61                  & 1032.21               \\ \cline{2-10}
\multicolumn{1}{|c|}{}                                                                           & Case 3.                  & 18                 & 189.19                & 36                  & 511.53                & 42                  & 624.28                & 55                  & 799.94                \\ \hline
\multicolumn{1}{|c|}{\multirow{3}{*}{\begin{tabular}[c]{@{}c@{}}$n=200$\\ $m=200$\end{tabular}}} & Case 1.                  & 29                 & 334.06                & 62                  & 741.88                & 74                  & 892.02                & 86                  & 985.20                \\ \cline{2-10}
\multicolumn{1}{|c|}{}                                                                           & Case 2.                  & 26                 & 298.12                & 52                  & 626.73                & 63                  & 762.45                & 74                  & 860.35                \\ \cline{2-10}
\multicolumn{1}{|c|}{}                                                                           & Case 3.                  & 25                 & 271.59                & 49                  & 545.83                & 59                  & 666.79                & 72                  & 792.59                \\ \hline
\multicolumn{1}{|c|}{\multirow{3}{*}{\begin{tabular}[c]{@{}c@{}}$n=300$\\ $m=100$\end{tabular}}} & Case 1.                  & 13                 & 210.13                & 45                  & 779.68                & 56                  & 958.04                & 64                  & 1020.85               \\ \cline{2-10}
\multicolumn{1}{|c|}{}                                                                           & Case 2.                  & 13                 & 210.74                & 38                  & 675.39                & 45                  & 800.78                & 53                  & 905.01                \\ \cline{2-10}
\multicolumn{1}{|c|}{}                                                                           & Case 3.                  & 17                 & 176.22                & 33                  & 449.02                & 41                  & 595.87                & 50                  & 734.57                \\ \hline
\end{tabular}
\caption{ Numerical results of mspACM for the minimax  problem \eqref{example-alt}. $n={\rm dim}\,x,\ m={\rm dim}\,y$; $\epsilon$-the relative error (i.e., $\|(x^k;y^k)\|_2/\|(x^0;y^0)\|_2$); $T$- number of iterations; $t$-CPU time.}\label{tab2}
\end{table}
\section{Some Concluding Remarks}\label{Sec4}
 \setcounter{equation}{0}
\quad \, Nonsmooth convex-concave minimax optimization problems are an important class of optimization problems with many applications. However, there are few  numerical algorithms for solving this type of problems when the smooth parts in the objective function are not bilinear. We developed a majorized semi-proximal alternating coordinate method (mspACM) for solving a nonsmooth convex-concave minimax problem of the form  (\ref{cminimax}). We demonstrated the global convergence of the algorithm mspACM under mild assumptions without requiring strong convexity-concavity condition and the linear rate of convergence under the locally metrical subregularity of the solution mapping. Preliminary numerical results have been reported, which shows the efficiency of the proposed mspACM method.

There are many interesting problems worth considering. In this paper, we only tested mspACM for two types of examples. How is the performance of  mspACM for other types of minimax optimization problems?  Presently we only considered convex-concave minimax problem of the form  (\ref{cminimax}); {\it i.e.}, $K(\cdot,\cdot)$ is required to be smooth convex-concave. How to construct efficient numerical algorithms for solving the nonsmooth minimax problem (\ref{cminimax}) when $K(\cdot,\cdot)$ is not a  smooth convex-concave function?\\
\medskip \noindent
{\bf Acknowledgments.} The authors are very grateful to Dr. Wenxing Zhang  for providing the important reference \cite{Valkonen2014}.

\end{document}